\documentclass[12pt]{amsart}
\usepackage{ifthen,amssymb}
\usepackage[all]{xypic}

\newcommand{\catD}{{D}}
\newcommand{\dbc}[1]{{D}^b(#1)}
\newcommand{\dpc}[1]{{D}^+(#1)}
\newcommand{\dmc}[1]{{D}^-(#1)}

\newcommand{\cdbc}[1]{{D}^b_c(#1)}

\newcommand{\fmf}[3]{{\Phi^{#1}_{{\scriptscriptstyle #2\!\rightarrow\! #3}}}}
\newcommand{\sk}[1]{{\cO_{#1}}}

\newcommand{\Hom}{{\operatorname{Hom}}}
\newcommand{\Ext}{{\operatorname{Ext}}}

\newcommand{\SHom}{{\mathcal{H}om}}
\newcommand{\SExt}{{\mathcal{E}xt}}

\newcommand{\lotimes}{{\,\stackrel{\mathbf L}{\otimes}\,}}

\newcommand{\Id}{{\operatorname{Id}}}

\DeclareMathOperator{\Spec}{{Spec}}

\DeclareMathOperator{\supp}{{supp}}

\newcommand{\cF}{{\mathcal F}}

\newcommand{\calH}{{\mathcal H}}
\newcommand{\cI}{{\mathcal I}}
\newcommand{\cJ}{{\mathcal J}}

\newcommand{\cK}{{\mathcal K}}
\newcommand{\cM}{{\mathcal M}}
\newcommand{\cN}{{\mathcal N}}
\newcommand{\cO}{{\mathcal O}}

\newcommand{\bR}{{\mathbf R}}
\newcommand{\bL}{{\mathbf L}}

\newcommand{\cplx}[1]{{{\mathcal #1}^{\scriptscriptstyle\bullet}}}

\newcommand{\dcplx}[1]{{{\mathcal #1}^{\scriptscriptstyle\bullet\vee}}}
\newcommand{\dSHom}[1]{{\mathcal{H}om_{#1}^{\scriptscriptstyle\bullet}}}

\newcommand{\marginnote}[1]{\ifthenelse{\isodd{\thepage}}{\normalmarginpar}
{\reversemarginpar}\marginpar{\fbox{\parbox{24mm}{\sloppy\footnotesize #1}}}}

\newcommand{\iso}{{\,\stackrel {\textstyle\sim}{\to}\,}}
\newcommand{\cono}{\operatorname{Cone}}

\DeclareMathSymbol{\functor}{\mathbin}{AMSa}{"20}

\newcommand{\scplx}[1]{{{\mathcal #1}^{\scriptscriptstyle\bullet\sharp}}}

\newcommand{\bcplx}[1]{{{\mathcal #1}_{Z_x}^{\scriptscriptstyle\bullet}}}

\newcommand{\Dcplx}[1]{{{\mathcal D}_{#1}^{\scriptscriptstyle\bullet}}}

 \newtheorem{thm}{Theorem}[section]
 \newtheorem*{thm*}{Theorem}
 \newtheorem{cor}[thm]{Corollary}
 \newtheorem{lem}[thm]{Lemma}
 \newtheorem{prop}[thm]{Proposition}
 
 \theoremstyle{definition}
 \newtheorem{defin}[thm]{Definition}
 \newenvironment{defn}{\begin{defin}}{\hfill\hspace{1pt}$\triangle$\end{defin}}
 \theoremstyle{remark}
 \newtheorem{rema}[thm]{Remark}
 \newtheorem{exe}[thm]{Example}
 
\newenvironment{rem}{\begin{rema}}{\hfill\hspace{1pt}$\triangle$\end{rema}}

\numberwithin{equation}{section}

\begin{document}
\title[RELATIVE INTEGRAL FUNCTORS AND SINGULAR  PARTNERS]{RELATIVE INTEGRAL FUNCTORS FOR SINGULAR FIBRATIONS AND SINGULAR  PARTNERS}
\author[D. Hern\'andez Ruip\'erez]{Daniel Hern\'andez Ruip\'erez}
\email{ruiperez@usal.es}
\author[A.C. L\'opez Mart\'{\i}n]{Ana Cristina L\'opez Mart\'{\i}n}
\email{anacris@usal.es}
\author[F. Sancho de Salas]{Fernando Sancho de Salas}
\email{fsancho@usal.es}
\address{Departamento de Matem\'aticas and Instituto Universitario de F\'{\i}sica Fundamental y Matem\'aticas
(IUFFYM), Universidad de Salamanca, Plaza
de la Merced 1-4, 37008 Salamanca, Spain}
\date{\today}
\thanks {Work supported by research projects MTM2006-04779 (MEC)
and SA001A07 (JCYL)} 
\begin{abstract}   We study relative integral functors for singular schemes and
characterise those which preserve boundness  and those which have integral right
adjoints. We prove that  a relative integral functor is an equivalence if and only if its
restriction to every fibre is an equivalence. This allows us to construct a non-trivial
auto-equivalence of the derived category of an arbitrary genus one fibration with no
conditions on either the base or the total space and getting rid of the usual
assumption of irreducibility of the fibres. We also extend to Cohen-Macaulay
schemes the criterion of Bondal and Orlov for an integral functor to be  fully faithful
in characteristic zero and give a different criterion which is valid in arbitrary
characteristic. Finally, we prove that for projective schemes both the Cohen-Macaulay and the Gorenstein conditions are invariant under Fourier-Mukai functors.
\end{abstract}
\maketitle

{\small \tableofcontents }

\section*{Introduction}

The relevance of derived categories and
Fourier-Mukai functors in birational geometry is nowadays well known \cite{Bri02,Kaw02b}. One of the most important problems in
this context is the minimal model problem. There are examples
proving that derived categories have a nice behavior under some
birational transformations as blow-ups, flips and flops.
Furthermore, Bondal and Orlov conjectured \cite{BO02} that each flip (resp.~flop) relating two smooth varieties $X$ and $X^+$ should induce a fully
faithful functor (resp. an equivalence) from the bounded derived category of coherent sheaves on $X^+$, $\cdbc{X^+}$,  to the corresponding bounded derived category
$\cdbc{X}$. Since varieties with singularities have to be allowed
in the minimal model programme, one of the main problems this
programme leads to is the study of derived categories for singular
projective varieties. However, not too much attention has been
paid to singular varieties in the literature on the topic. The reason may be
that many of the fundamental results rely deeply on smoothness.

In dimension three, the difficulties arising when one deals with
Gorenstein schemes can be circumvent using a smoothing approach
\cite{Chen02}. Threefold Gorenstein terminal singularities are isolated
hypersurface singularities (\cite{KoMo98}) and a
hypersurface singularity can be considered as a special fiber of a smooth
fourfold. Thus, one can get some information about the
derived categories in dimension 3 from the derived categories of
the corresponding smooth fourfold. The results in \cite{Chen02} can be
generalized  to $Q$-Gorenstein teminal threefolds (see \cite{AC05, Kaw02b,
Kaw02a}). In this situacion, the associated Gorenstein stack allows to
reduce the problem to the Gorenstein case. Nevertheless, since in higher
dimension there are not local models for quotient singularities
and these seem to be more rigid, this smoothing idea does not work
and new results are needed.

Following a completely different path, we started in \cite{HLS05} the study of derived categories and fully
faithful integral functors for schemes with Gorenstein
singularities or fibered in schemes of that kind, and generalised to that situation the characterisation of fully faithfulness originally proved by Bondal and Orlov \cite{BO95} in the smooth case.

The aim of this paper is twofold. On the one hand, we give a very
general result that characterizes when relative integral functors
are equivalences and that allow to reduce the problem to the
absolute setting. This gives in particular a construction of a
non-trivial invertible integral functor for a general genus one
fibration. Notice that we adopt here a slightly different definition of a genus  one fibration; actually we mean a flat Gorenstein morphism whose fibres are curves of arithmetic genus 1 and have trivial dualising sheaf; the last condition is a consequence of the others  when the fibres are reduced.
In the case of dimensions two or three, we are then
allowing all types of Kodaira fibre, thus getting rid of the usual
assumption of irreducible fibres.

Moreover, for fully faithfulness
we are now able to consider more general singularities, namely
Cohen-Macaulay schemes, both in the absolute and in the relative
case. The case of positive characteristic, never studied before,
is considered as well.  On the second hand, we also tackle the
question of what geometric information of the scheme can be
recovered from its derived category. For projective schemes, we
give an affirmative answer for both the Cohen-Macaulay and the
Gorenstein properties.

The paper is organised as follows:

Section \ref{basicSec} is a repository of formulas used throughout
the work and collected here for the reader's convenience.

In
Section \ref{intfuncSec} we study integral functors, which we define
 directly in the relative situation, and determine conditions for them
 to map bounded complexes to bounded complexes and to have right adjoints
 with the same property. The study requires the  notions of complexes of
 relative finite projective dimension and of relative finite homological
 dimension, which are equivalent in the case of projective morphisms as
 we proved in \cite{HLS05}. Among the results in Section \ref{intfuncSec}
 we can mention Proposition \ref{p:equivfhd}, where we prove that if the
 integral functor defined by a kernel $\cplx K\in \cdbc{X\times Y}$ induces
 an equivalence between the bounded derived categories of $X$ and $Y$,
 then the kernel has to be of finite homological dimension over both $X$
 and $Y$ and the same happens to the kernel $\bR \dSHom{\cO_{X\times Y}}(\cplx K, \pi_Y^!\cO_Y)$
 of the right adjoint.

The main result of this Section is however Proposition
\ref{p:relative} that proves, in great generality, that a relative
integral functor is fully faithful (or an equivalence) if and only
if the absolute integral functors induced on the fibers are fully
faithful (or equivalences). This result could be useful for the
study of relative Fourier-Mukai transforms between the derived
categories of the total spaces of two very general fibrations. The
reason is that it allows to pass from a relative to an absolute
situation where more things are known. For
instance in the absolute setting, there is an important class of
well-known Fourier-Mukai functors: the twist functors along
spherical objects that were firstly introduced by Seidel and
Thomas in \cite{SeTh01}. As a direct corollary of this
proposition, we construct a non-trivial auto-equivalence of the
derived category of an arbitrary genus one fibration. The result
is the following:

\begin{thm*}[Proposition \ref{p:poincequiv}] Let $S$ be an algebraic scheme,
$X\to S$ a genus one fibration, and $\cI_\Delta$ the ideal sheaf
of the relative diagonal. The relative integral functor
$$\fmf{\cI_\Delta}{X}{X}\colon \cdbc{X}\to \cdbc{X}$$ is an equivalence of categories.
\end{thm*}

This extends the result and the proof given in \cite[Prop.~2.7]{HLS05} in two directions. First the characteristic of the
base field is arbitrary. Second, all kind of possible fibres of a
genus one fibration are allowed; in particular, for dimensions two
or three they are all the Kodaira fibres. Moreover, no condition is imposed on either the base $S$ or the total space $X$ of the fibration.

When the fibration has only integral fibres, Proposition \ref{p:poincequiv} gives a short
proof of the invertibility of the usual elliptic integral functor. We then expect that Proposition \ref{p:poincequiv} could be a useful tool for the study of the moduli spaces of relatively semistable sheaves on $X\to S$ with respect to a suitable relative ample divisor following \cite{BBHM02,BBH08,BuKr04,HePl05}, and for the study of the derived category of $X$ generalising \cite{BuKr05}.

In Section \ref{ffSec} we give criteria to characterise fully
faithful integral functors. We first consider the case of
characteristic zero, and extend to Cohen-Macaulay schemes the
characterisation of fully faithful integral functors given by
Bondal and Orlov. The characterisation result is the following.

\begin{thm*}[Theorem \ref{1:ffcritCM}]  Let $X$ and $Y$ be proper schemes over
an algebraically closed field of characteristic zero, and let
$\cplx{K}$ be an object in $\cdbc{X\times Y}$ of finite
homological dimension over both $X$ and $Y$. Assume also that $X$
is projective, Cohen-Macaulay and integral. Then the functor
$\fmf{\cplx{K}}{X}{Y}\colon \cdbc{X}\to \cdbc{Y}$ is fully
faithful if and only if the kernel $\cplx{K}$ is strongly simple
over $X$.
\end{thm*}

This theorem also generalises the extension to varieties with
Gorenstein singularities given in \cite[Thm.~1.22]{HLS05}. The
new theorem is also stronger in the Gorenstein case, because we do
not need to assume any longer that $Y$ is projective and
Gorenstein.

As in the Gorenstein case, strong simplicity  in the Cohen-Macaulay case
 (Definition \ref{1:strgsplcplxCM}) is defined in terms of locally complete
intersection zero cycles instead of the structure sheaves of the
closed points.  In the smooth case, our definition is weaker that
the usual one given by Bondal and Orlov, and then  Theorem
\ref{1:ffcritCM} improves the characterization of fully
faithfulness of Bondal and Orlov. It should be noticed that we
give here a slightly different notion of strong simplicity than
the one given in \cite{HLS05} for Gorenstein varieties. Since both
characterise the fully faithfulness of the corresponding integral
functors, they are equivalent in the Gorenstein case.

We also consider the case of positive characteristic which is different because
Theorem \ref{1:ffcritCM} fails in that situation. We have modified the notion of
strong simplicity to a new one, which also characterises fully faithfulness in
arbitrary characteristic. The precise statement is:

\begin{thm*}[Theorem \ref{1:ffcritCMp}]  Let $X$ and $Y$ be proper schemes over
 an algebraically closed field of arbitrary characteristic, and
 let $\cplx{K}$ be an object
in $\cdbc{X\times Y}$ of finite homological dimension over both
$X$ and $Y$. Assume also that $X$ is connected, equidimensional,
projective and Cohen-Macaulay. Then the functor
$\fmf{\cplx{K}}{X}{Y}\colon \cdbc{X}\to \cdbc{Y}$ is fully
faithful if and only if the kernel $\cplx{K}$ has the following
properties:
\begin{enumerate}
\item For every  closed point $x\in X$ there is a l.c.i.~zero cycle
$Z_x$ supported on $x$ such that
$$
\Hom^i_{\catD(Y)}(\fmf{\cplx{K}}{X}{Y}(\cO_{Z_{x_1}}),\fmf{\cplx{K}}{X}{Y}(\cO_{x_2}))=0
$$
unless $x_1= x_2$ and $0\leq i\leq \dim X$.
\item There exists a closed point $x$ such that at least one of the following conditions is fulfilled:
\begin{enumerate} \item $\Hom^0_{\catD(Y)}(\fmf{\cplx{K}}{X}{Y}(\cO_X),
\fmf{\cplx{K}}{X}{Y}(\sk{x}))\simeq k$.
\item $\Hom^0_{\catD(Y)}(\fmf{\cplx{K}}{X}{Y}(\cO_{Z_x}),
\fmf{\cplx{K}}{X}{Y}(\sk{x}))\simeq k$ for any l.c.i.~zero cycle
$Z_x$ supported on $x$.
\item $\dim_k \Hom^0_{\catD(Y)}(\fmf{\cplx{K}}{X}{Y}(\cO_{Z_x}),
\fmf{\cplx{K}}{X}{Y}(\cO_{Z_x}))\leq l(\cO_{Z_x})$ for any
l.c.i.~zero cycle $Z_x$ supported on $x$, where $l(\cO_{Z_x})$ is
the length of $\cO_{Z_x}$.
\end{enumerate}
\end{enumerate}

\end{thm*}

Due to Proposition \ref{p:relative} and the properties of relative
integral functors proved in Section \ref{intfuncSec}, the
extension of the above criteria to the relative setting is
straightforward (see Theorem \ref{t:relative}).



As an application we give a  different proof of Proposition
\ref{p:poincequiv} which doesn't need to use the work of Seidel
and Thomas about spherical objects.

The last Section \ref{intpartnersSec} is devoted to the study of Fourier-Mukai
partners of a given proper scheme $X$, that is, proper schemes $Y$  with equivalent
coherent bounded category (i.e, $D$-equivalent to $X$), and such that the equivalence
is given by an integral functor.  In the projective smooth case, the second condition
is automatically fulfilled, due to  Orlov's representation theorem  \cite{Or97}. However,
the validity of Orlov's theorem for singular varieties is still unknown; then, in principle,
two $D$-equivalent  singular varieties might not be Fourier-Mukai partners.

It is known that smooth Fourier-Mukai partners share many geometrical properties.
In the same vein, we prove in Section \ref{intpartnersSec} that singular Fourier-Mukai
partners also have many geometrical properties in common. Our main result in this direction is the following.

\begin{thm*}[Theorem \ref{t:intpartnersCM}]
Let $X$ be a projective equidimensional Cohen-Macaulay scheme and $Y$ a projective
Fourier-Mukai partner of $X$. Then one has
\begin{enumerate}
\item If $Y$ is reduced, then $Y$ is equidimensional and $\dim Y=\dim X$.
\item If $Y$ is equidimensional and $\dim Y=\dim X$, then $Y$ is Cohen-Macaulay. Moreover, if $X$ is
Gorenstein, then $Y$ is Gorenstein as well.
\end{enumerate}
\end{thm*}

\subsubsection*{Conventions}
In this paper, scheme means separated  scheme of
finite type over an algebraically closed field $k$. By a
Gorenstein or a Cohen-Macaulay morphism, we understand a flat
morphism of schemes whose fibres are respectively Gorenstein or
Cohen-Macaulay. For any scheme $X$ we denote by $\catD(X)$ the
derived category of complexes of $\cO_X$-modules with
quasi-coherent cohomology sheaves. This is the essential image of
the derived category of quasi-coherent sheaves in the derived
category $\catD(\mathfrak{Mod}(X))$ of all $\cO_X$-modules \cite[Cor.~5.5]{BoNee93}. Analogously $\dpc{X}$, $\dmc{X}$
and $\dbc{X}$ will denote the derived categories of complexes
which are respectively bounded below, bounded above and bounded on
both sides, and have quasi-coherent cohomology sheaves. The
subscript $c$ will refer to the corresponding subcategories of
complexes with coherent cohomology sheaves.

\subsection*{Acknowledgements} We would like to thank to the authors of the forthcoming book \cite{BBH08} for sharing with us their notes and to the anonymous referees for comments and suggestions which helped us to improve the manuscript.
Ana Cristina L\'opez Mart\'{\i}n would like also to thank Miles Reid for useful comments and for his warm welcome in Warwick.

\section{Some basic formulas in derived category}\label{basicSec}


We recall here some basic formulas which will be used in the rest
of the paper.

If $X$ is a scheme, there is  a functorial isomorphism (in the
derived category)
\begin{equation}
\label{tens1}
 \bR \dSHom{\cO_X} (\cplx{F},\bR \dSHom{\cO_X} (\cplx{E},\cplx{H})) \iso
 \bR \dSHom{\cO_X} (\cplx{F}\lotimes \cplx{E}, \cplx{H})
 \end{equation}
  where
$\cplx{F}$,  $\cplx{E}$ and  $\cplx{H}$ are in $\catD(\mathfrak{Mod}(X))$ \cite[Thm.~A]{Spal98}. One
also has a functorial isomorphism in $\catD(\mathfrak{Mod}(X))$
\begin{equation}
\label{tens2} \bR \dSHom{\cO_X} (\cplx{F},\cplx{E}) \lotimes
\cplx{H} \iso \bR \dSHom{\cO_X} (\cplx{F},\cplx{E}\lotimes
\cplx{H})
\end{equation}
when either $\cplx F$ or $\cplx H$ has finite
homological dimension. When $\cplx F$ is bounded above with coherent cohomology sheaves, $\cplx E$ is bounded below and $\cplx H$ has finite homological dimension, the formula is standard (cf.~\cite[Prop.~II.5.14]{Hart66}). Since we have not found a reference for the unbounded case, we give a simple proof here:
Let $\cplx{E}\to\cplx{I}$ be an injective resolution \cite[Thm.~4.5]{Spal98} and  $\cplx{P}\to\cplx{H}$ a flat resolution \cite[Prop.~5.6]{Spal98}. One has morphisms of complexes
$$
\dSHom{\cO_X} (\cplx{F},\cplx{I}) \otimes \cplx{P}\to  \dSHom{\cO_X} (\cplx{F},\cplx{I} \otimes \cplx{P})\to
\dSHom{\cO_X} (\cplx{F},\cplx{J} )
$$
where $\cplx{J}$ is an injective resolution of $\cplx{I} \otimes \cplx{P}$. This proves the existence of a morphism
$$
\bR \dSHom{\cO_X} (\cplx{F},\cplx{E}) \lotimes
\cplx{H} \to \bR \dSHom{\cO_X} (\cplx{F},\cplx{E}\lotimes \cplx{H})\,.
$$
We now prove that this is an isomorphism if either $\cplx F$ or $\cplx H$ has finite homological dimension. This is a local question, so we may assume that $X=\Spec A$ is affine and that either $\cplx{F}$ or
$\cplx{P}$ is a bounded complex of free $A$-modules of finite rank. In both cases we have an isomorphism of complexes
$$
\dSHom{A} (\cplx{F},\cplx{I})\otimes \cplx{P} \iso \dSHom{A} (\cplx{F},\cplx{I}\otimes \cplx{P} )\,.
$$
In the first case, this proves directly the isomorphism \eqref{tens2}. In the second case, we have only to take into account that since $\cplx{P}$ is a bounded complex of free modules of finite rank, the complex $\cplx{I}\otimes \cplx{P}$ is injective.

If $f\colon X\to Y$ is a morphism of schemes, the direct and inverse images are defined for unbounded complexes, $\bR f_\ast\colon \catD(\mathfrak{Mod}(X))\to \catD(\mathfrak{Mod}(Y))$, $\bL f^\ast\colon \catD(\mathfrak{Mod}(Y))\to \catD(\mathfrak{Mod}(X))$  and the latter is a right adjoint to the former \cite[Thm.~B]{Spal98}. They induce morphisms  $\bR f_\ast\colon \catD(X)\to \catD(Y)$, $\bL f^\ast\colon \catD(Y)\to \catD(X)$ which are adjoint of each other as well. One has a ``projection formula''
$$
\bR f_\ast(\cplx F\lotimes\bL f^\ast\cplx G) \simeq \bR f_\ast\cplx F\lotimes\cplx G
$$
for $\cplx F$ in $\catD(X)$ and $\cplx G$ in $\catD(Y)$ \cite[Prop.~5.3]{Nee96}.

Let us consider a cartesian diagram of morphisms of algebraic varieties
$$
\xymatrix{X\times_Y Z  \ar[r]^(.7){\tilde g} \ar[d]^{\tilde f}& X\ar[d]^f \\
Z  \ar[r]^g & Y }
$$
Then for any complex $\cplx G$ of $\cO_X$-modules there is a natural
morphism
$$
\bL g^\ast\bR f_\ast \cplx G\to \bR \tilde f_\ast \bL {\tilde
g}^\ast \cplx G\,.
$$
Moreover, if $\cplx G$ is in $\catD(X)$ and either $f$ or $g$
is flat, the above morphism is an isomorphism. This is the so-called ``base-change formula'' in the derived category. The flat base-change formula, i.e., when $g$ is flat, is well-known \cite[Prop.~II.5.12]{Hart66}. If  $g$ is arbitrary and $f$ is flat, the formula is proven in \cite[Appendix A]{BBH08}. In this paper we only need the following very simple case.
\begin{prop}\label{p:basechange}
Let us consider a diagram
$$
\xymatrix{X_y \ar@{^(->}[r]^{j_{X_y}} \ar[d]^{f_y}& X\ar[d]^f \\
\{y\}  \ar@{^(->}[r]^{j_y} & Y }
$$
where $f$ is a flat morphism of schemes, $y\in Y$ is a closed point and $X_y=f^{-1}(y)$ is the fibre. For every object $\cplx G$ in $\catD(X)$, there is a base-change isomorphism
$$
\bL j_y^\ast\bR f_\ast \cplx G\simeq \bR  f_{y\ast} \bL j_{X_y}^\ast \cplx G\,.
$$
in the derived category.
\end{prop}
\begin{proof}
It is enough to prove that the induced morphism 
$$
j_{y\ast}(\bL j_y^\ast\bR f_\ast \cplx G)\to j_{y\ast}(\bR  f_{y\ast} \bL j_{X_y}^\ast \cplx G)
$$
 is an isomorphism in $\catD(Y)$. By the projection formula, the first member is isomorphic to $\cO_y\lotimes  \bR f_\ast \cplx G$. The second member is isomorphic to
$\bR f_\ast j_{X_y \ast} \bL j_{X_y}^\ast \cplx G$, which is isomorphic to $\bR f_\ast (\cO_{X_y}\lotimes\cplx G)
$
again by the projection formula. Moreover $\cO_{X_y}\simeq \bL f^\ast\cO_y$ because $f$ is flat, and then $\bR f_\ast (\cO_{X_y}\lotimes\cplx G)\simeq \cO_y\lotimes \bR f_\ast\cplx G$ by the projection formula for $f$.
\end{proof}

Let $f\colon X\to Y$ be a proper morphism of schemes. The relative
Grothendieck duality states the existence of a functorial
isomorphism in the derived category
\begin{equation}
\label{1:duality} \bR\dSHom{\cO_Y}(\bR f_{\ast}\cplx{F},\cplx{G})
\simeq \bR f_\ast \bR\dSHom{\cO_{X}} (\cplx{F}, f^!\cplx{G})\,.
\end{equation} for $\cplx G$ in $\catD(Y)$ and $\cplx F$ in
$\catD(X)$ (see for instance \cite{Nee96}). By applying the
derived functor of the global section functor, we obtain the
\emph{global duality formula}
\begin{equation} \label{1:adj1}
\Hom_{\catD(Y)}(\bR f_{\ast}\cplx{F},\cplx{G}) \simeq
\Hom_{\catD(X)} (\cplx{F}, f^!\cplx{G})\,.
\end{equation}
In other words, the direct image $\bR f_{\ast}\colon \catD(X)\to
\catD(Y)$ has a right adjoint $ f^!\colon \catD(Y) \to \catD(X)$. If $g\colon Y \to Z$ is another
proper morphism, there is a natural functor isomorphism $(g\circ f)^!\simeq f^!\circ g^!$.

We shall call the objet $f^!\cO_Y$ the
\emph{relative dualizing complex} of $X$ over $Y$. When $Y$ is a
point, we also write $\Dcplx{X}$ instead of $f^!\cO_Y$.

When $f$ is flat, then it is Cohen-Macaulay if and only if the relative dualizing complex $f^!\cO_Y$ is isomorphic to a single sheaf $\omega_{X/Y}$ placed at degree $-n$ (where $n$ is the relative dimension of $f$). Moreover $\omega_{X/Y}$ is a line bundle if and only if $f$ is a Gorenstein morphism.

We finish this part by recalling some properties of Grothendieck duality which we shall use in this paper.

Firstly, if $f\colon X\to Y$ is a finite morphism, then $f_\ast f^!\cplx G
\simeq\bR\dSHom{\cO_Y}(f_\ast \cO_X,\cplx G)$.
Secondly, Grothendieck duality is compatible with flat base-change,
that is,  if $g\colon
Z\to Y$ is a flat morphism and $f_Z\colon Z\times_Y X\to Z$ is the
induced morphism, then
$f_Z^! \cO_Z  \simeq g_X^\ast f^!\cO_Y$ where $g_X\colon Z\times_Y
X\to X$ is the projection. In particular, the formation of $f^!\cO_Y$ is compatible with open immersions $U\subseteq Y$.
Finally, there is a natural map $\bL f^\ast 
\cplx{G}\lotimes f^!\cO_Y\to f^!\cplx{G}$; in some cases, it is an isomorphism. 
One of those cases is when $\cplx{G}$ has finite homological dimension; to see this one has to prove that for any complex $\cplx F$ in $\catD(X)$ the induced map $\Hom_{\catD(X)}(\cplx F, \bL f^\ast 
\cplx{G}\lotimes f^!\cO_Y) \to \Hom_{\catD(X)}(\cplx F, f^!\cplx{G})$ is bijective, and for this it is enough to prove that $\bR f_\ast \bR\dSHom{\cO_X}(\cplx F, \bL f^\ast 
\cplx{G}\lotimes f^!\cO_Y) \to \bR f_\ast\bR\dSHom{\cO_X}(\cplx F, f^!\cplx{G})$ is an isomorphism in $\catD(Y)$. 
The first member is isomorphic to $\bR\dSHom{\cO_Y}(\bR f_\ast \cplx F,\cO_Y)\lotimes\cplx G$, by Equation \eqref{tens2}, 
the projection formula and the duality isomorphism \eqref{1:duality}.
Moreover, $\bR\dSHom{\cO_Y}(\bR f_\ast \cplx F,\cO_Y)\lotimes\cplx G 
\simeq \bR f_\ast\bR\dSHom{\cO_X}(\cplx F, f^!\cplx{G})$ by Equation \eqref{tens2} and relative duality  \eqref{1:duality}.
Other case where $\bL f^\ast 
\cplx{G}\lotimes f^!\cO_Y\to f^!\cplx{G}$ is an isomorphism is when 
$f$ is a regular closed immersion; in this case one has to prove that the induced morphism $f_\ast(\bL f^\ast  \cplx{G}\lotimes f^!\cO_Y)\to f_\ast(f^!\cplx{G})$ is an isomorphism. By the projection formula the first member is isomorphic to $ \cplx{G}\lotimes f_\ast(f^!\cO_Y)$.
Since  $f_\ast f^!\cO_Y
\simeq\bR\dSHom{\cO_Y}(f_\ast \cO_X,\cO_Y)$ and $f_\ast f^!\cplx G
\simeq\bR\dSHom{\cO_Y}(f_\ast \cO_X,\cplx G)$, the result follows again from Equation \eqref{tens2}, which can be applied because $f_\ast \cO_X$ is of finite homological dimension.

%

\section{Relative Integral functors for singular schemes} \label{intfuncSec}
\subsection{Boundedness conditions}

Let $S$ be a scheme and let $X\to S$ and $Y\to S$ be proper
morphisms. We denote by $\pi_X$ and $\pi_Y$ the projections of the
fibre product $X\times_SY$ onto its factors.

Let $\cplx{K}$ be an object in $\cdbc{X\times_SY}$. The relative
integral functor defined by $\cplx{K}$ is the functor $\fmf{\cplx
K}{X}{Y} \colon \catD (X) \to \catD (Y)$ given by $$\fmf{\cplx
K}{X}{Y}(\cplx{F})=\bR \pi_{Y\ast}( \bL\pi_X^\ast\cplx{F}\lotimes
\cplx{K})\, .$$ and it maps $\dmc{X}$ to $\dmc{Y}$.

By adjunction between the direct and inverse images and by
relative Grothendieck duality, $\fmf{\cplx K}{X}{Y}$ has a  right
adjoint $H\colon \catD(Y)\to \catD(X)$ given by
\begin{equation}\label{e:adjduality}
H(\cplx{G})=\bR \pi_{X,\ast}( \bR\dSHom{\cO_{X\times Y}}
(\cplx{K}, \pi_Y^!\cplx{G}))\,,
\end{equation}
 which maps $\dpc{Y}$ to $\dpc{X}$. We shall now study when either $\fmf{\cplx
K}{X}{Y}$ or $H$ take bounded complexes to bounded complexes.

\begin{defn} Let $f\colon Z \to T$ be a morphism of
schemes. An object $\cplx E$ in $\cdbc{Z}$ is said to be of
\emph{finite homological dimension} (resp. of  \emph{finite
projective dimension}) \emph{over $T$}, if $\cplx E\lotimes \bL
f^\ast \cplx G$ (resp. $\bR\dSHom{\cO_{X}} (\cplx{E},
f^!\cplx{G})$) is bounded for any $\cplx{G}$ in $\cdbc{T}$.
\end{defn}

\begin{rem} The absolute notion of finite homological dimension corresponds to being ``of finite homological dimension over $Z$ with respect to the identity'' \cite[Lemma 1.2]{HLS05}, rather than to being ``of finite homological dimension over $\Spec k$''; actually, any object $\cplx E$ in $\cdbc{Z}$ is of finite homological dimension over $\Spec k$. The usual notion of morphism of finite homological dimension is equivalent to saying that $\cO_Z$ is ``of finite homological dimension over $T$''.
\end{rem}

For projective morphisms one has the following result (see
\cite{HLS05}):

\begin{prop}\label{fpdchar2}  Let $f\colon Z\to T$ be a projective morphism and $\cO(1)$ a relatively
very ample line bundle. Let $\cplx{E}$ be an object of $\cdbc{Z}$
and let us set $\cplx E(r)=\cplx E \otimes \cO(r)$. The following
conditions are equivalent:
\begin{enumerate}
\item   $\cplx E$ is of finite projective dimension over $T$.
\item $\bR f_\ast
(\cplx E(r))$ is of finite homological dimension (i.e. a perfect complex) for every
integer $r$.
\item $\cplx E$ is of finite homological dimension over $T$.
\end{enumerate}
Thus, if $f$ is locally projective, $\cplx E$ is of finite
projective dimension over $T$ if and only if it is of  finite
homological dimension over $T$.\qed\end{prop}

\begin{cor} Let $f\colon Z\to T$ be a locally projective morphism. The functor $f^!$ sends $\dbc{T}$ to $\dbc{Z}$ if and only if $f$ is a morphism of finite homological dimension.
\qed
\end{cor}

We also state here the following lemma whose proof can be
found in \cite{HLS05}.
\begin{lem} \label{l:beilinson}
Let $f\colon Z \to T$ be a projective morphism.
\begin{enumerate}
\item Let $\cplx E$ be an object of  $\dmc{Z}$. Then $\cplx E=0$ (resp. is an object of $\dbc{Z}$) if
and only if $\bR f_\ast (\cplx E(r))=0$ (resp. is an object of
$\dbc{T}$) for every integer $r$.
\item Let $g\colon \cplx E\to\cplx F$ be a morphism in $\dmc{Z}$. Then $g$ is an isomorphism if and only if the
induced morphism $\bR f_\ast (\cplx E(r))\to \bR f_\ast (\cplx
F(r))$ is an isomorphism in $\dmc{T}$ for every integer $r$.
\end{enumerate}
\qed\end{lem}

The arguments used in the proof of the above Lemma also show that an object $\cplx E$ of $\catD(Z)$ has coherent cohomology sheaves if and only if $\bR f_\ast (\cplx E(r))$ has coherent cohomology sheaves for every integer $r$. One then obtains the following result.
\begin{prop} If $f\colon Z \to T$ is a locally projective morphism, the functor $f^!$ sends $\catD_c(T)$ to $\catD_c(Z)$.
\end{prop}
\begin{proof} Let $\cplx G$ be an object of $\catD_c(T)$. Since the formation of $f^!\cplx G$ is compatible with open immersions, we can assume that $f$ is projective. In this case,  Grothendieck duality gives $\bR f_\ast (f^!\cplx G(r))\simeq \bR\dSHom{\cO_Z}(\bR f_{\ast}\cO_Z(-r),\cplx{G})$, and the latter complex has coherent  cohomology sheaves. We finish by Lemma \ref{l:beilinson}.
\end{proof}

\begin{prop} \label{p:kernelfhd}
Assume that $X\to S$ is locally projective  and let $\cplx{K}$ be
an object in $\cdbc{X\times_S Y}$. The functor
$\fmf{\cplx{K}}{X}{Y}$ maps $\cdbc{X}$ to $\cdbc{Y}$ if and only
if $\cplx{K}$ has finite homological dimension over $X$.
\end{prop}
\begin{proof} Assume that  $\fmf{\cplx{K}}{X}{Y}$ maps $\cdbc{X}$ to $\cdbc{Y}$. We have to prove that
$\bL \pi_X^\ast\cplx F\lotimes\cplx K$ is bounded for any bounded
complex  $\cplx F$. We may assume that $X\to S$ is projective, and
then $\pi_Y$ is also projective. By Lemma \ref{l:beilinson}, it suffices
to show that $\bR\pi_{Y\ast}[\bL \pi_X^\ast\cplx F\lotimes\cplx
K\otimes \pi_X^\ast\cO (r)]$ is bounded for any $r$. This is
immediate from the equality $\bR\pi_{Y\ast}[\bL \pi_X^\ast\cplx
F\lotimes\cplx K\otimes \pi_X^\ast\cO (r)]= \fmf{\cplx{K}}{X}{Y}
(\cplx F (r))$. The converse is clear.
\end{proof}

\begin{prop}
\label{dualGore} Let $f\colon Z\to T$ be a locally projective
morphism of  schemes and $\cplx{E}$  an object of $\cdbc{Z}$ of
finite homological dimension over $T$. One has
$$
\bR \dSHom{\cO_{Z}} (\cplx{E}, f^!\cO_T)\lotimes \bL f^\ast \cplx
G\simeq \bR \dSHom{\cO_{Z}} (\cplx{E}, f^!\cplx{G})
$$
for $\cplx{G}$ in $\cdbc{T}$. In particular, $\bR \dSHom{\cO_{Z}}
(\cplx{E}, f^!\cO_T)$ is also  of finite homological dimension
over $T$.
\end{prop}

\begin{proof} One has natural morphisms
\begin{equation}
\bR\dSHom{\cO_Z}(\cplx{E},  f^!\cO_T)\lotimes \bL
f^\ast\cplx{G}\to
\bR\dSHom{\cO_Z}(\cplx{E}, \bL f^\ast\cplx{G}\lotimes f^!\cO_T) \\
\to \bR\dSHom{\cO_Z}(\cplx{E}, f^! \cplx{G})\,. \label{eqiso}
\end{equation}
 We have to prove that the composition is an isomorphism. This is
a local question on $T$, so that we can assume that  $f$ is
projective.

By Lemma \ref{l:beilinson} we have to prove that the induced morphism
\begin{equation}
\bR f_{\ast}(\bR\dSHom{\cO_Z}(\cplx{E},f^!\cO_T)\lotimes \bL
f^\ast\cplx{G} \otimes \cO (r)) \to \bR
f_{\ast}(\bR\dSHom{\cO_Z}(\cplx{E}, f^! \cplx{G})\otimes \cO (r))
\end{equation}
is an isomorphism for any integer $r$. By Grothendieck duality,
this is equivalent to proving that the induced morphism
\begin{equation}
\bR\dSHom{\cO_T}(\bR f_\ast \cplx{E}(-r),\cO_T)\lotimes \cplx{G}
\to \bR\dSHom{\cO_T}(\bR f_\ast \cplx{E}(-r), \cplx{G})\,,
\label{eqiso2}
\end{equation}
is an isomorphism. This follows from \eqref{tens2} because $\bR
f_\ast \cplx{E}(-r)$ is of finite homological dimension by
Proposition \ref{fpdchar2}. Finally, by Proposition
\ref{fpdchar2}, $\bR \dSHom{\cO_{Z}} (\cplx{E}, f^!\cO_T)$ is also
of finite homological dimension over $T$.
\end{proof}

\begin{prop} \label{p:AdjGo} Assume that $X\to S$ is locally projective and let
$\cplx{K}$ be an object in $\cdbc{X\times_S Y}$ of finite
homological dimension over both $X$ and $Y$. The functor
$\fmf{\bR\dSHom{\cO_{X\times_S Y}}(\cplx{K}, \pi_Y^!\cO_Y) }{Y}{X}
\colon \cdbc{Y}\to\cdbc{X}$ is a right adjoint  to
$\fmf{\cplx{K}}{X}{Y}\colon \cdbc{X}\to\cdbc{Y}$.
\end{prop}
\begin{proof} Since $\pi_Y$ is a locally projective morphism,
Proposition \ref{dualGore} gives the result.
\end{proof}

When the integral functor $\fmf{\cplx K}{X}{Y}$ is an equivalence
between the bounded  categories, Proposition \ref{p:AdjGo} can be
refined as follows:
\begin{prop} \label{p:equivfhd} Let $X\to S$ and $Y\to S$ be projective morphisms and let $\cplx K$
be a kernel in $\cdbc{X\times_S Y}$ such that $\fmf{\cplx
K}{X}{Y}$ induces an equivalence
$$
\fmf{\cplx K}{X}{Y}\colon \cdbc{X}\simeq \cdbc{Y}\,.
$$
Then one has:
\begin{enumerate}
\item $\cplx K$ is finite of homological dimension over both $X$ and $Y$.
\item The right adjoint to $\fmf{\cplx K}{X}{Y}\colon \cdbc{X}\simeq \cdbc{Y}$ is an integral
functor, and its kernel $\bR\dSHom{\cO_{X\times_S Y}}(\cplx
K,\pi_Y^!\cO_Y)$  is also of finite homological dimension
over both $X$ and $Y$.
\end{enumerate}
\end{prop}
\begin{proof} First, $\cplx K$ is of finite homological dimension over $X$ by Proposition
\ref{p:kernelfhd}. We still don't know if $\cplx K$ is of finite
homological dimension over  $Y$, so we cannot apply Proposition
\ref{p:AdjGo}. However, we can proceed as follows: The functor
$\fmf{\cplx K}{X}{Y}\colon \catD(X)\to \catD(Y)$ over the whole
derived category has a right adjoint $H$ (cf.~Equation
\ref{e:adjduality}). Since $\fmf{\cplx K}{X}{Y}\colon
\cdbc{X}\simeq \cdbc{Y}$ is an equivalence, it has a right adjoint
$\bar H \colon \cdbc{Y}\to \cdbc{X}$. Let us prove that  $H$ and
$\bar H$ coincide over $\cdbc{X}$. For any $\cplx F\in \cdbc {X}$,
$\cplx G\in \cdbc{Y}$ one has that
$$\Hom_{\catD (X)}( \cplx F , \bar H(\cplx G)) \simeq \Hom_{\catD (Y)}
(\fmf{\cplx K}{X}{Y}(\cplx F),\cplx G) \simeq \Hom_{\catD (X)} (
\cplx F ,H(\cplx G))$$ Hence, there is a morphism $\eta\colon \bar
H(\cplx G)\to H(\cplx G)$ in $\dpc{X}$, such that the induced
morphism
$$
\Hom^i_{\catD(X)}(\cplx F,\bar H(\cplx G))\to
\Hom^i_{\catD(X)}(\cplx F,H(\cplx G))
$$
is an isomorphism for any $\cplx F\in \cdbc{X}$ and every integer
$i$. Let $\cplx C\in \dpc{X}$ be the cone of $\eta$. Then
$\Hom^i_{\catD(X)}(\cplx F, \cplx C)=0$ for any $\cplx F\in
\cdbc{X}$ and any $i$. Taking $\cplx F = \cO_X(r)$ for $r$ big
enough, one concludes that $\cplx C=0$, so that $\bar H \simeq H$.

Now, for any $\cplx G\in\cdbc{Y}$ one has that
$$
 \bar H(\cplx G\otimes\cO(r)) \simeq \bR\pi_{X\ast}[ \bR\dSHom{\cO_{X\times_S Y}} (\cplx K, \pi_Y^!\cplx
G)\otimes \pi_Y^\ast \cO(r)]
$$
 is bounded for any integer $r$. By Lemma \ref{l:beilinson},
$\bR\dSHom{\cO_{X\times_S Y}}(\cplx K, \pi_Y^!\cplx G)$ is bounded
and then $\cplx K$ has
finite homological dimension over $Y$. It follows from Proposition
\ref{p:AdjGo} that $\bar H$ is an integral functor of kernel
$\bR\dSHom{\cO_{X\times_S Y}}(\cplx K,\pi_Y^! \cO_Y)$. Since $\bar H$ is
an equivalence, its kernel is also of finite homological
dimension over both $X$ and $Y$.
\end{proof}

\begin{prop} \label{dsharp}
Let $f\colon Z\to T$ be a locally projective morphism and
$\cplx{E}$ an object of $\cdbc{Z}$. If $\cplx{E}$ has finite
homological dimension over $T$ then
 the natural morphism $$ \cplx{E}\to \bR \dSHom{\cO_{Z}}
(\bR \dSHom{\cO_{Z}} (\cplx{E}, f^!\cO_T), f^!\cO_T)$$ is an
isomorphism.
\end{prop}
\begin{proof} The problem is local on $T$, so we can assume that $f$ is projective.
Since $\bR \dSHom{\cO_{Z}} (\cplx{E}, f^!\cO_T)$ is of finite
homological dimension over $T$   by Proposition \ref{dualGore},
the object $\bR \dSHom{\cO_{Z}} (\bR \dSHom{\cO_{Z}} (\cplx{E},
f^!\cO_T), f^!\cO_T)$ has bounded cohomology. Thus, by Lemma
\ref{l:beilinson}, to prove that
$$
\cplx E \to \bR \dSHom{\cO_{Z}} (\bR \dSHom{\cO_{Z}} (\cplx{E},
f^!\cO_T), f^!\cO_T)
$$ is an isomorphism we have to prove that the induced morphism
$$
\bR f_\ast (\cplx E(r)) \to \bR f_\ast(\bR \dSHom{\cO_{Z}} (\bR
\dSHom{\cO_{Z}} (\cplx{E}, f^!\cO_T), f^!\cO_T)(r))
$$
is an isomorphism for every integer $r$. Grothendieck duality
implies that this morphism is the natural morphism $\bR f_\ast
(\cplx E(r)) \to (\bR f_\ast (\cplx E(r)))^{\vee\vee}$ (where for
a complex $\cplx F$ we write $\dcplx F=\bR \dSHom{\cO_T}(\cplx
F,\cO_T)$ for the derived dual), which is an isomorphism because
$\bR f_\ast (\cplx E(r)) $ is of finite homological dimension by
Proposition \ref{fpdchar2}.
\end{proof}

If $X$ is a proper scheme over a field, we denote by $\scplx E$
the ``dual in the dualising complex'', that is,
\begin{equation}\label{e:sharpdef}
\scplx E=\bR\dSHom{\cO_X}(\cplx E,
\Dcplx{X})\,.
\end{equation}
 If $X$ is projective, we can apply Proposition
\ref{dsharp} to the projection $f\colon X\to \Spec k$. Since any
bounded complex $\cplx E$ has finite homological dimension over
$\Spec k$ one has:

\begin{cor}  \label{dsharpequiv} Let $X$ be a projective scheme
over a field. The contravariant functor $\sharp$ induces an anti-equivalence of
categories $\sharp\colon \cdbc{X} \to \cdbc{X}$. \qed\end{cor}

We finish this subsection with the following property about
commutation of integral functors with the functor $\sharp$ (cf.~Equation \ref{e:sharpdef}), which will be used in
Theorems \ref{1:ffcritCM} and \ref{t:intpartnersCM}.

\begin{lem} \label{sharpfmf}
Let $X$ and $Y$ be projective schemes over a field and
$\cplx{K}\in \cdbc{X\times Y}$ a kernel of finite homological
dimension over $X$. Then one has
$$
(\fmf{\cplx K}{X}{Y}(\cplx{E}))^\sharp \simeq
\fmf{\bR \dSHom{\cO_{X\times Y}}(\cplx K,\pi_X^!\cO_X)}{X}{Y}(\scplx{E})\,.
$$
for any object $\cplx E$ in $\cdbc{X}$.

\end{lem}
\begin{proof} On the one hand, we have
\begin{align*}
(\fmf{\cplx K}{X}{Y}(\cplx{E}))^\sharp &\simeq \bR\dSHom{\cO_Y}
(\bR \pi_{Y\ast}(\pi_X^\ast \cplx{E}\lotimes \cplx K),\Dcplx{Y}) \simeq
\bR \pi_{Y\ast}\bR\dSHom{\cO_{X\times Y}}(\pi_X^\ast \cplx{E}\lotimes
\cplx K, \pi_Y^!\Dcplx{Y})\\
& \simeq \bR \pi_{Y\ast}\bR\dSHom{\cO_{X\times Y}} (\cplx K, \bR\dSHom{\cO_{X\times Y}}(\pi_X^\ast \cplx{E},
\pi_Y^!\Dcplx{Y}))\,.
\end{align*}
On the other hand
\begin{align*}
\fmf{\bR \dSHom{\cO_{X\times Y}}(\cplx K,\pi_X^!\cO_X)}{X}{Y}(\scplx{E}) &\simeq \bR \pi_{Y\ast}(
\pi_X^\ast \bR\dSHom{\cO_X}(\cplx{E},\Dcplx{X})\lotimes \bR\dSHom{\cO_{X\times Y}}
(\cplx K,\pi_X^!\cO_X))
\\ &\simeq \bR\pi_{Y\ast} \bR\dSHom{\cO_{X\times Y}}(\cplx K, \pi_X^!\bR\dSHom{\cO_X}(\cplx{E},\Dcplx{X}))
\end{align*}
where the last isomorphism is by Proposition \ref{dualGore}. Then, it is enough to see that
$\bR\dSHom{\cO_{X\times Y}}(\pi_X^\ast \cplx{E},
\pi_Y^! \Dcplx{Y})\simeq \pi_X^!\bR\dSHom{\cO_X}(\cplx{E},\Dcplx{X})$. Since $\pi_Y^! \Dcplx{Y}\simeq
\Dcplx{X\times Y}\simeq \pi_X^! \Dcplx{X}$, this follows
from the isomorphisms
\begin{align*}
 \Hom_{\dbc{X\times Y}} (\cplx F, \bR\dSHom{\cO_{X\times Y}}(\pi_X^\ast \cplx{E}, \Dcplx{X\times Y}))&\simeq
 \Hom_{\dbc{X\times Y}} (\cplx F\lotimes
 \pi_X^\ast \cplx{E}, \Dcplx{X\times Y})\\ &\simeq
\Hom_{\dbc{X}}(\bR \pi_{X\ast}(\cplx F\lotimes \pi_X^\ast \cplx{E}), \Dcplx{X})\\ &\simeq
\Hom_{\dbc{X}}(\bR \pi_{X\ast}\cplx F\lotimes \cplx{E}, \Dcplx{X}) \\ &\simeq
\Hom(\bR \pi_{X\ast}\cplx F,\bR\dSHom{\cO_X}(\cplx{E}, \Dcplx{X}))\\ & \simeq \Hom_{\dbc{X\times Y}}(\cplx F,
\pi_X^!\bR\dSHom{\cO_X}(\cplx{E},\Dcplx{X}))
\end{align*}
which hold for any object $\cplx F$ in $\cdbc{X\times Y}$.
\end{proof}

\subsection{Restriction to fibres: a criterion for equivalence}

In this subsection we prove that, in a very general situation, to
see that a relative integral functor is fully faithful (or an
equivalence) it is enough to prove that its restriction to each
fibre is fully faithful (or an equivalence). Using this result and
the theorem of Seidel and Thomas that proves that any twist
functor along a spherical object is an equivalence of categories,
we construct then a non-trivial auto-equivalence of the bounded
derived category of an arbitrary genus one fibration.

The base field is here an algebraically closed field of arbitrary
characteristic. Let $p\colon X\to S$ and $q\colon Y\to S$ be
proper and \emph{flat} morphisms. Let $\cplx K$ be an object in
$\cdbc{X\times_S Y}$ and $\Phi=\fmf{\cplx{K}}{X}{Y}$. For any
closed point $s\in S$ we write $X_s=p^{-1}(s)$, $Y_s=q^{-1}(s)$,
and denote by $\Phi_s\colon \dmc{X_s} \to \dmc{Y_s}$ the integral
functor defined by $\cplx{K}_s=\bL j_s^\ast \cplx{K}$, with
$j_s\colon X_s\times Y_s\hookrightarrow X\times_S Y$ the natural
embedding.

When the kernel $\cplx{K}\in \cdbc{X\times_SY}$ is of finite
homological dimension over $X$, the functor $\Phi$ maps $\cdbc{X}$ into $\cdbc{Y}$ and $\cplx{K}_s\in \cdbc{X_s\times Y_s}$ for any $s\in S$. Morever,
since $q$ is flat, $\cplx{K}_s$ is of finite homological dimension
over $X_s$.

From the base-change formula (Proposition \ref{p:basechange}, see also \cite{BBH08}) we obtain that
\begin{equation}\label{e:basechange}\bL
j_s^\ast\Phi (\cplx{F})\simeq\Phi_s(\bL j_s^\ast\cplx{F})
\end{equation}
for every $\cplx{F}\in \catD(X )$, where $j_s\colon
X_s\hookrightarrow X$ and $j_s\colon Y_s\hookrightarrow Y$ are the
natural embeddings.  In this situation, base change formula also
gives that
\begin{equation}\label{e:directimage}
j_{s\ast}\Phi_s(\cplx{G})\simeq\Phi ( j_{s\ast}\cplx{G})
\end{equation}
for every $\cplx{G}\in \catD(X_s)$.

\begin{lem} \label{nonbound} Let $\Psi\colon \dmc{X}\to \dmc{Y}$ be an integral
functor whose kernel is an object of  $\dmc{X\times Y}$. For any integer $i$ there exists
$r_0$ such that $\calH^i (\Psi (\cplx{G})) \simeq \calH^i
(\Psi(\sigma_{\geq r} \cplx{G}) )$ for every $r\leq r_0$ and any
$\cplx G \in \dmc{X}$, where $\sigma_{\geq r} \cplx{G} $ is the
truncation (cf.~\cite[I\S 7]{Hart66}).
\end{lem}
\begin{proof} There exist an integer $m$ such that if $\cplx F$ is an object of $\dmc{X}$ and
$\calH^j(\cplx F)=0$ for $j\geq s$, then $\calH^j(\Psi(\cplx
F))=0$ for $j\geq s+m$. Let us take
$r\leq i-1-m$ and let $\alpha_r\colon \cplx{G}\to \sigma_{\geq
r}(\cplx G)$ be the natural morphism. Then
$\calH^j(\cono(\alpha_r))=0$ for $j\geq r$, so that
$\calH^j(\Psi(\cono(\alpha_r)))=0$ for $j\geq r+m$. In particular
$\calH^{i-1}(\Psi(\cono(\alpha_r)))=\calH^{i}(\Phi^-(\cono(\alpha_r)))=0$,
and we finish by taking cohomology on  the exact triangle
$$
\Psi(\cplx G) \to \Psi(\sigma_{\geq r}(\cplx G)) \to
\Psi(\cono(\alpha_r))\to \Psi(\cplx G)[1]\,.
$$
\end{proof}

\begin{prop}  \label{p:relative}
Assume that $X\to S$ is locally projective and let $\cplx{K}$ be
an object in $\cdbc{X\times_SY}$ of finite homological dimension
over both $X$ and $Y$. The relative integral functor
$\Phi=\fmf{\cplx{K}}{X}{Y}\colon \cdbc{X}\to \cdbc{Y}$ is fully
faithful (resp.~an equivalence) if and only if $\Phi_s\colon
\cdbc{X_s}\to \cdbc{Y_s}$ is fully faithful (resp.~an equivalence)
for every closed point $s\in S$.
\end{prop}
\begin{proof}
By Proposition \ref{p:AdjGo}, the integral functor
$H=\fmf{\bR\dSHom{X\times Y}(\cplx{K},\pi_Y^!\cO_Y)}{Y}{X} \colon
\cdbc{X}\to \cdbc{Y}$ is a right adjoint to $\Phi$. We can now
proceed as in the proof of \cite[Thm.~2.4]{HLS05}, which we explain here in some more detail.

 If $\Phi$ is
fully faithful the unit morphism $\Id \to H\circ \Phi$ is an
isomorphism. Then, given a closed point $s\in S$ and $\cplx G\in
\cdbc{X_s}$, one has an isomorphism $j_{s\ast}\cplx G\to (H\circ
\Phi)(j_{s\ast}\cplx G)$.  Since $(H\circ \Phi)(j_{s\ast}\cplx
G)\simeq j_{s\ast} (H_s\circ \Phi_s) (\cplx G)$ by
\eqref{e:directimage} and $j_s$ is a closed immersion, the unit
morphism $\cplx G\to (H_s\circ \Phi_s)(\cplx G)$ is an
isomorphism; this proves that $\Phi_s$ is fully faithful.

Now assume that $\Phi_s$ is fully faithful for any closed point
$s\in S$.
To prove that $\Phi$ is fully faithful we have to see that the unit morphism $\eta\colon\Id\to H\circ\Phi$ is
an isomorphism. For each $\cplx F\in \cdbc{X}$ we have an exact
triangle
$$
\cplx F\xrightarrow{\eta(\cplx F)} (H\circ \Phi)(\cplx F)\to \cono(\eta (\cplx F))\to \cplx F[1]\,.
$$ 
Let us now fix a closed point $s\in S$ . By Equation \eqref{e:basechange}, we have an exact triangle
$$
\bL j_s^\ast \cplx F\to (H_s\circ \Phi_s)(\bL j_s^\ast \cplx F)\to \bL
j_s^\ast \cono(\eta (\cplx F))\to \bL j_s^\ast \cplx F[1]\,.
$$
By Lemma \ref{nonbound}, for every integer $i$ there exists $r$ small enough such that
$\calH^i(\bL j_s^\ast \cplx F)\simeq \calH^i(\sigma_{\geq r}\bL
j_s^\ast \cplx F)$ and $\calH^i((H_s\circ \Phi_s)(\bL j_s^\ast
\cplx F))\simeq \calH^i((H_s\circ \Phi_s)(\sigma_{\geq r}\bL
j_s^\ast \cplx F))$. Since $\sigma_{\geq r}\bL j_s^\ast \cplx F$ is
a bounded complex, and $\eta_s\colon \Id \to
H_s\circ \Phi_s$ is an isomorphism because $\Phi_s$ is fully faithful, one has that $\calH^i(\sigma_{\geq r}\bL j_s^\ast \cplx F)\simeq
\calH^i[(H_s\circ \Phi_s)(\sigma_{\geq r}\bL j_s^\ast \cplx F)]$. Thus $\bL j_s^\ast \cplx F\to (H_s\circ \Phi_s)(\bL j_s^\ast \cplx F)$ induces isomorphisms between all the cohomology sheaves, so that it is an isomorphism, and then $\bL
j_s^\ast \cono(\eta (\cplx F))=0$. Since this holds for every closed point $s\in S$, we finish  by  \cite[Lemma 2.3]{HLS05}.


A similar argument gives the statement about equivalence.
\end{proof}

Related results concerning Azumaya smooth varieties have been proved by Kuznetsov in  \cite[Prop.~2.44 and Thm.~2.46]{Kuz06}.  To apply his results to our situation, $X$, $Y$ and $S$ have to be smooth though the flatness conditions on $p$ and $q$ can be removed.

\subsubsection{An auto-equivalence of the derived category of a genus one
fibration}\label{ss:autoequiv} Let $p\colon X\to S$ be a genus one
fibration, that is, a projective Gorenstein morphism whose fibres
are curves with arithmetic genus $\dim H^1(X_s,\cO_{X_s})=1$ and
have trivial dualising sheaf. No further assumptions on $S$ or $X$
are made here.

When the fibres are reduced of arithmetic genus one, then the
condition on the dualising sheaf is always fulfilled. However,
since nonreduced curves can also appear as degenerated fibres for
such a genus one fibration, and for these curves the dualizing
sheaf need not to be trivial (see for instance \cite{CFHR99}), one
needs to assume it.

There are some cases where the structure of the singular fibers is
known: For smooth elliptic surfaces over the complex numbers, the
classification was given by Kodaria  \cite{kod} and for smooth
elliptic threefolds over a base field of characteristic different
from 2 and 3, they were classified by Miranda \cite{Mir83}. In
both cases, the possible singular fibres are plane curves of the
same type, the so-called Kodaira fibres. Nevertheless, in a genus
one fibration non-plane curves could appear as degenerated
fibres.

Since we are not putting any restriction on the characteristic of
the base field or on the dimension of $X$, our genus one fibrations may have singular fibres other than the Kodaira fibres.

We consider the commutative diagram
$$\xymatrix{ & X\times_S X \ar[ld]_{\pi_1}\ar[rd]^{\pi_2} \ar[dd]^{\rho}&  \\
X \ar[rd]^{p} & & X\ar[ld]_{p} \\
& S & }$$ and the relative integral functor
$$
\Phi=\fmf{\cI_\Delta}{X}{X}\colon \cdbc{X}\to  \cdbc{X}\,,
$$
with kernel the ideal sheaf $ \cI_\Delta$ of the relative diagonal
immersion $\delta\colon X \hookrightarrow X\times_S X$.

\begin{prop} \label{p:poincequiv}
The relative integral functor
$$\Phi=\fmf{\cI_\Delta}{X}{X}\colon \cdbc{X}\to \cdbc{X}
$$  defined by the ideal sheaf of the relative diagonal is an equivalence
of categories.
\end{prop}

\begin{proof} To prove that $\cI_\Delta$
is of finite homological dimension over both factors,  it is
enough to see that it is of finite homological dimension over the
first factor because of its symmetry. By the exact sequence
$$
0\to \cI_\Delta\to \cO_{X\times_S X}\to \delta_\ast \cO_X\to 0 \,,
$$
it suffices to see that $\delta_\ast \cO_X$ has finite homological
dimension over the first factor. We have then to prove that for
any $\cplx{N}\in \dbc{X}$, the complex $\delta_\ast \cO_X\otimes
\pi_1^\ast\cplx{N}$ is bounded and this follows from the
projection formula for $\delta$.

For every closed point $s\in S$, the absolute functor
$\Phi_s=\fmf{\cI_{\Delta_s}}{X_s}{X_s}$ is equal to
$T_{\mathcal{O}_{X_s}}[-1]$ where $T_{\mathcal{O}_{X_s}}$ denotes
the twist functor along the object $\mathcal{O}_{X_s}$. Since
$X_s$ is a genus one projective curve with trivial dualizing
sheaf, $\mathcal{O}_{X_s}$ is a spherical object. Thus,
$T_{\mathcal{O}_{X_s}}$ is an equivalence of categories by
\cite{SeTh01} and we conclude by Proposition \ref{p:relative}.
\end{proof}

Similar results has been obtained by Burban and Kreussler in \cite{BuKr06}. They proved a  version of Proposition \ref{p:poincequiv} in the case when the base field is of characteristic zero, $S$ and $X$ are
reduced, $X$ is connected, the fibration $p\colon X\to S$ has only
integral fibres and it has a section taking values in the smooth locus of $X$. None of these assumptions have been made in this section, though we assume by technical reasons that $S$ (and then $X$) is separated.

Notice that, since we are working in the relative setting, the
integral functor $\fmf{\cI_\Delta}{X}{X}$ whose kernel is the
ideal sheaf of the relative diagonal is not easily described as a
twist functor. Actually, the twist functor $T_{\cO_X}[-1]$ is the
integral functor whose kernel is the ideal sheaf of the absolute
diagonal immersion $X\hookrightarrow X\times X$. Even the latter
functor may fail to be an equivalence because in general
$\mathcal{O}_X$ is not spherical.

We shall give an alternative proof of this result, without using
the work of Seidel and Thomas, in Subsection \ref{ellfibss}.

\section{Fully faithfulness criteria for integral functors} \label{ffSec}
\subsection{Strongly simple objects and spanning classes}

Recall that if $X$ is a Cohen-Macaulay scheme, for every point $x$
there exist a zero cycle $Z_x$ supported on $x$ defined locally by
a regular sequence (cf.~\cite[Lemma 1.9]{HLS05}); we refer to such
cycles as  to \emph{locally complete intersection} or
\emph{l.c.i.}~cycles. If $Z_x\hookrightarrow X$ is
 a l.c.i.~cycle, by the Koszul complex theory, the structure
sheaf $\cO_{Z_x}$ has \emph{finite homological dimension} as an
$\cO_X$-module.

In order to fix some notation, for any zero-cycle $Z_x$ of $X$ and
any scheme $S$, we shall denote by $j_{Z_x}$ the immersion
$Z_x\times S \hookrightarrow X\times S$.

For further use, we gather here two equivalent characterisations of objects of the derived
category of an equidimensional scheme (that is, a scheme with all its irreducible components of the same dimension) defined by single sheaves supported on a closed subscheme.  The statements are a slight generalisation of  \cite[Props.~1.15, 1.16]{HLS05} and the proofs are completely analogous.

\begin{prop} Let $X$, $Y$ be equidimensional Cohen-Macaulay schemes, $j\colon Y\hookrightarrow X$ a closed immersion
of codimension $d$, and $\cplx{K}$ an object of $\cdbc X$. Assume that
\begin{enumerate}
\item If $x\in X-Y$ is a closed point, then $\bL j_{Z_x}^\ast \cplx{K}=0$ for some
l.c.i.~zero cycle $Z_x$ supported on $x$.
 \item If $x\in Y$ is a closed point, then $L_i
j_{Z_x}^\ast \cplx{K}=0$ for some l.c.i.~zero cycle $Z_x$ supported on $x$ when
either $i<0$ or $i>d$.
\end{enumerate}
Then there is a sheaf $\cK$ on $X$ whose topological support is
contained in $Y$ and  such that $\cplx{K}\simeq\cK$ in $\cdbc X$. Moreover, if $\cplx K$ is nonzero, then the topological support of $\cK$ is a union of irreducible components of $Y$.
\label{1:support} \qed \end{prop}

\begin{prop}  \label{1:support2} Let $X$, $Y$ be equidimensional Cohen-Macaulay schemes of dimensions $m$ and $n$ respectively, $j\colon Y\hookrightarrow X$ a closed immersion, and $\cplx{K}$ an
object of $\cdbc X$. Assume that for any closed point $x\in X$
there is a l.c.i.~zero cycle $Z_x$ supported on $x$ such that
$$
 \Hom^i_{\catD(X)}(\cO_{Z_x},\cplx{K})= 0\,,
$$
unless $x\in Y$ and $m\le i\le n$. Then there is a sheaf $\cK$ on $X$ whose topological support is contained in $Y$ such that $\cplx{K}\simeq\cK$ in $\cdbc X$.
Moreover, if $\cplx K$ is nonzero, then the topological support of $\cK$ is a union of irreducible components of $Y$.
\qed\end{prop}

In this paper we only use the first part of the statement, namely the fact that $\cplx{K}\simeq\cK$ and its topological support is contained in $Y$.

Let $X$ and $Y$ be proper schemes. Assume that $X$ is Cohen-Macaulay.
In this situation, the notion of strong simplicity is the following.

\begin{defn} An object  $\cplx{K}$ in $\cdbc{X\times Y}$ is \emph{strongly simple} over $X$ if it satisfies the
following conditions:
\begin{enumerate}
\item For every  closed point $x\in X$ there is a l.c.i.~zero cycle
$Z_x$ supported on $x$ such that
$$
\Hom^i_{\catD(Y)}(\fmf{\cplx{K}}{X}{Y}(\cO_{Z_{x_1}}),\fmf{\cplx{K}}{X}{Y}(\cO_{x_2}))=0
$$
unless $x_1= x_2$ and $0\leq i\leq \dim X$.
\item$\Hom^0_{\catD(Y)}(\fmf{\cplx{K}}{X}{Y}(\sk{x}),
\fmf{\cplx{K}}{X}{Y}(\sk{x}))= k$ for every closed point $x\in X$.
\end{enumerate}
\label{1:strgsplcplxCM}
\end{defn}

The last condition can be written as $\Hom^0_{\catD(Y)}(\bL
j_{x}^\ast\cplx{K}, \bL j_{x}^\ast\cplx{K})= k$, because the
restriction $\bL j_x^\ast\cplx{K}$ of $\cplx{K}$ to the fibre $j_x\colon Y\simeq \{x\}\times Y\hookrightarrow
X\times Y$ can
also be  computed  as $\fmf{\cplx{K}}{X}{Y}(\sk{x})$. Then, the objects $\bL j_{x}^\ast\cplx{K}$ are simple and we may think of $\cplx K$ as a parametrisation of simple objects of $\cdbc{Y}$; this is the reason why one calls $\cplx K$ a strongly simple object over $X$.

\begin{rem}
When $X$ and $Y$ are smooth, our definition is weaker than the
usual one given by Bondal and Orlov (see \cite{BBH08}). As a
consequence of a result of  Bondal and Orlov
\cite[Thm.~1.1]{BO95}, and of Theorem \ref{1:ffcritCM}, our
definition is equivalent to this one in the smooth case. In
\cite{HLS05} we gave another notion of strongly simple objects for
Gorenstein schemes. Again, Theorem \ref{1:ffcritCM} and the
corresponding statement \cite[Theorem 1.22]{HLS05} for the
Gorenstein  case, prove that the two notions are equivalent in
that situation.
\end{rem}

\subsubsection{Spanning classes}
As in the Gorenstein
case, the derived category $\cdbc{X}$ has a natural spanning class
when $X$ is a Cohen-Macaulay scheme.

\begin{lem} \label{l:spanning} If $X$ is a Cohen-Macaulay scheme, then
the set
$$\Omega=\{ \cO_{Z_x} \text{ for all closed points $x\in X$ and all
l.c.i.~zero cycles $Z_x$ supported on $x$} \}$$ is a spanning
class for $\cdbc{X}$.
\end{lem}
\begin{proof} Take a non-zero object $\cplx E$ in $\cdbc{X}$. On the one hand,
for every l.c.i.~zero cycle as above  there is a spectral sequence
$E_2^{p,q}=\Ext_{\cO_X}^p(\calH^{-q}(\cplx E),\cO_{Z_x})$
converging to $E_2^{p+q}=\Hom^{p+q}_{\catD (X)}(\cplx
E,\cO_{Z_x})$. If $q_0$ is the maximum of  the $q$'s such that
$\calH^{q}(\cplx E)\neq 0$ and $x$ is a point of the support of
$\calH^{q}(\cplx E)$, then $E_2^{0,-q_0}\neq 0$ for every
l.c.i.~zero cycle $Z_x$ and any non zero-element there survives to
infinity. Then $\Hom^{-q_0}_{\catD (X)}(\cplx E,\cO_{Z_x})\neq 0$.

On the other hand, by Proposition \ref{1:support2} with
$Y=\emptyset$, if $\Hom^{i}_{\catD (X)}(\cO_{Z_x},\cplx E)= 0$ for
every $i$ and every $Z_x$, then $\cplx{E}=0$.
\end{proof}

\subsection{A criterion in characteristic zero}

We now give the criterion for an integral functor between derived
categories of Cohen-Macaulay proper schemes to be fully faithful.


\begin{thm}\label{1:ffcritCM}  Let $X$ and $Y$ be proper schemes over
an algebraically closed field of characteristic zero, and let
$\cplx{K}$ be an object in $\cdbc{X\times Y}$ of finite
homological dimension over both $X$ and $Y$. Assume also that $X$
is projective, Cohen-Macaulay and integral. Then the functor
$\fmf{\cplx{K}}{X}{Y}\colon \cdbc{X}\to \cdbc{Y}$ is fully
faithful if and only if the kernel $\cplx{K}$ is strongly simple
over $X$.
\end{thm}
\begin{proof}  If the functor $\fmf{\cplx{K}}{X}{Y}$ is fully faithful, then $\cplx{K}$ is
strongly simple over $X$.

Let us prove the converse. Before starting, we fix some notation:
we denote by $\pi_i$ the projections of $X\times X$ onto its
factors and $\Phi =\fmf{\cplx{K}}{X}{Y}$.

By Proposition \ref{p:AdjGo}, the integral functor
$H=\fmf{\bR\dSHom{X\times Y}(\cplx{K},\pi_Y^!\cO_Y)}{Y}{X}$ is a right
adjoint to $\Phi$. By \cite[Prop.~1.18]{HLS05} it suffices to show
that $H\circ \Phi $ is fully faithful. We know that the
composition of integral functors is again an integral functor (cf.~\cite{HLS05}), and
then $H\circ \Phi \simeq \fmf{\cplx{M}}{X}{X}$, with ${\cplx
M}\in \cdbc{X\times X}$. 

The strategy of the proof is similar to that of \cite[Thm.~1.22]{HLS05}. We are going to prove that
$\cplx M\simeq \delta_\ast \cN$ where $\delta \colon X\hookrightarrow X\times X$ is the diagonal immersion
and $\cN$ is a line bundle on $X$; then $\fmf{\cplx{M}}{X}{X}$ consist of twisting by $ \cN$ which is an
equivalence of categories, in particular fully faithful.

\medskip\noindent
a) \emph{$\bL j_x^\ast \cplx M$ is  a single sheaf topologically
supported on $x$. Thus $\cplx{M}$ is a single sheaf $\cM$ topologically supported on
the diagonal and flat over $X$ by the first projection by \cite[Lemma 4.3]{Bri99}.}

\medskip

Let us fix a closed point $(x_1,x_2)\in X\times X$ and consider
the l.c.i.~zero cycle $Z_{x_1}$ of the first condition of the
definition of strongly simple object.  One has
$$
\Hom^i_{\catD(X)}(\cO_{Z_{x_1}}, \fmf{\cplx{M}}{X}{X}(\cO_{x_2}))
\simeq \Hom^i_{\catD(Y)}(\Phi(\cO_{Z_{x_1}}),\Phi(\cO_{x_2}))\,,
$$
which is zero unless $x_1=x_2$ and $0\le i\le m$ because $\cplx K$
is strongly simple.  Applying Proposition \ref{1:support2} to  the
immersion $\{x_2\}\hookrightarrow X$ we have that
$\fmf{\cplx{M}}{X}{X}(\cO_{x_2})$ reduces to a
coherent sheaf topologically supported on $x_2$. Since
$Lj_x^\ast \cM\simeq \fmf{\cplx{M}}{X}{X}(\cO_{x})$, one has that
$\cplx M$ is a sheaf $\cM$ whose support is contained in the
diagonal and $\pi_{1\ast}\cM$ is locally free, where $\pi_1\colon
X\times X\to X$ is the  projection onto the first factor. The rank
of $\pi_{1\ast}\cM$ can not be zero, by condition (2) of strongly
simple. Hence, the topological support of $\cplx M$ is exactly the
diagonal.

\medskip\noindent
b) \emph{$\cM$ is schematically supported on the diagonal, that
is, $\cM = \delta_\ast \cN$ for a coherent sheaf $\cN$ on $X$;
moreover $\cN$ is a line bundle.}

\medskip
Let us denote by $\bar\delta\colon W \hookrightarrow X\times
X$ the schematic support of $\cM$ so that $\cM = \bar \delta_\ast
\cN$ for a coherent sheaf $\cN$ on $W$.  Since $\cM$ is topologically supported on the diagonal,
the diagonal embedding
$\delta$ factors through a closed immersion $\tau\colon
X\hookrightarrow W$ which topologically is a homeomorphism.

Since $\cM$ is flat over $X$ by $\pi_1$, $\cN$ is flat over $X$ by the composition
$\bar \pi_1=\pi_1\circ \bar\delta\colon W \to X$. Moreover $\bar\pi_1$ is a finite morphism,
so that  $\pi_{1\ast}\cM\simeq \bar\pi_{1\ast} \cN$ is locally free.

Now, as in the proof of \cite[Thm.~1.22]{HLS05}, to conclude  is enough to see that
the rank of the locally free sheaf $\pi_{1\ast}\cM$ is one.

One has that $\Hom^0(\cO_x,\fmf{\cM}{X}{X}(\cO_x))\simeq
\Hom^0(\Phi(\cO_x),\Phi(\cO_x))\simeq k$. Hence there is an
injective morphism $\cO_x\to \fmf{\cM}{X}{X}(\cO_x)\simeq
j_x^\ast \cM$. It suffices to show that this is an isomorphism for at
least one closed point $x$. If $C_x$ is the cokernel, we have to
see that $\Hom^0(\cO_x, C_x)=0$. Arguing like in \cite{Bri99} it
is enough to show that $\Hom^1(\cO_x,\cO_x)\to
\Hom^1(\fmf{\cM}{X}{X}(\cO_x),\fmf{\cM}{X}{X}(\cO_x))$ is
injective.

Let us denote $\widetilde\cM =\bR\dSHom{X\times X}(\cM,\pi_1^!\cO_X)$, which
is again a sheaf supported topologically on the diagonal and flat
over $X$ by the first projection, because
$$\pi_{1\ast}\bR\dSHom{X\times Y}(\cM,\pi_1^!\cO_X)\simeq
\bR\dSHom{X\times Y}(\pi_{1\ast}\cM, \cO_X) $$ and $\pi_{1\ast}\cM$ is
locally free.

One has that $\cO_x^\sharp\simeq \cO_x$ by relative duality for the closed immersion $\{x\}\hookrightarrow X$, so that
$\bL j_x^\ast \widetilde \cM \simeq (\bL j_x^\ast \cM)^\sharp$ by  Lemma \ref{sharpfmf}. Then,
$$
\Hom^0(\bL
j_x^\ast \widetilde\cM,\cO_x)\simeq   \Hom^0((\bL
j_x^\ast \cM)^\sharp,\cO_x^\sharp)\simeq \Hom^0(\cO_x,\bL j_x^\ast \cM)\simeq
k
$$
because the functor $\sharp$ is an anti-equivalence of categories by Corollary \ref{dsharpequiv}. Hence,
$\bL j_x^\ast \widetilde\cM$ is the
sheaf of a zero cycle supported on $x$. We can then apply \cite[Lemmas 5.2, 5.3]{Bri99} to $\widetilde\cM$ to
conclude that there exists a point $x$ such that
\begin{equation}\label{e:kodaira}
\Hom^1(\cO_x,\cO_x)\to
\Hom^1(\fmf{\widetilde\cM}{X}{X}(\cO_x),\fmf{\widetilde\cM}{X}{X}(\cO_x))\simeq
\Hom^1(\fmf{\cM}{X}{X}(\cO_x)^\sharp,\fmf{\cM}{X}{X}(\cO_x)^\sharp)
\end{equation}
is injective.

Now, again by Corollary  \ref{dsharpequiv}, we have a commutative diagram
$$
\xymatrix{
\Hom^1(\cO_x,\cO_x) \ar[d]^{\sharp}_{\simeq} \ar[r] & \Hom^1(\fmf{\cM}{X}{X}(\cO_x),\fmf{\cM}{X}{X}(\cO_x))
\ar[d]^{\sharp}_{\simeq}\\
\Hom^1(\cO_x,\cO_x) \ar[r]  & \Hom^1(\fmf{\cM}{X}{X}(\cO_x)^\sharp,\fmf{\cM}{X}{X}(\cO_x)^\sharp)\,.
}
$$
Since the bottom arrow is injective, the top arrow is injective as well and we conclude.
\end{proof}

The condition that the characteristics of the base field $k$ is zero is only used to prove that there exists a point $x$ such that the morphism \eqref{e:kodaira} is injective \cite[Lemmas 5.2, 5.3]{Bri99}. This is proving by showing that \eqref{e:kodaira} is the Kodaira-Spencer map for the family $\widetilde\cM$.  Moreover, the sheaves $\widetilde\cM_x$ define points of a Hilbert scheme and the above Kodaira-Spencer map is the composition of the tangent map to the map $x\mapsto \widetilde\cM_x$, and the Kodaira-Spencer map for the universal family. The latter is an isomorphism and since $k$ has characteristics zero and $x\mapsto \widetilde\cM_x$ is injective on closed points, its tangent map is injective at some point. It is this last statement what fails to be true in positive characteristics.


\subsection{A criterion in arbitrary characteristic}

As we showed in \cite{HLS05}, Theorem \ref{1:ffcritCM} is no
longer true in positive characteristic even in the smooth case. We
reproduce here the counterexample given there. Let $X$ be a smooth
projective scheme of dimension $m$ over a field $k$ of
characteristic $p>0$, and $F\colon X \to X^{(p)}$ the
\emph{relative Frobenius morphism} \cite[3.1]{Illu96}, which is
topologically a homeomorphism.  Let $\Gamma\hookrightarrow X\times
X^{(p)}$ be the graph of $F$, whose associated integral functor is
the direct image $F_\ast \colon \cdbc{X} \to \cdbc{X^{(p)}}$.
Since $F_\ast(\cO_x)\simeq \cO_{F(x)}$, one easily sees that
$\Gamma$ is strongly simple over $X$. However, $F_\ast(\cO_{X})$
is a locally free $\cO_{X^{(p)}}$-module of rank $p^m$
\cite[3.2]{Illu96}, so that
$\Hom^0_{\catD(X^{(p)})}(F_\ast(\cO_X), \cO_{F(x)})\simeq k^{p^m}$
whereas $\Hom^0_{\catD(X)}(\cO_{X},\cO_x)\simeq k$; thus $F_\ast$
is not fully faithful.

Then, in  arbitrary characteristic we need another characterisation of those kernels which give rise to fully faithful integral functors. The right notion is the following

\begin{defn} \label{d:ortho}
An object $\cplx K$ of $\cdbc{X\times Y}$ satisfies the Cohen-Macaulay orthonormality conditions over $X$ if  it has the following
properties:
\begin{enumerate}
\item For every  closed point $x\in X$ there is a l.c.i.~zero cycle
$Z_x$ supported on $x$ such that
$$
\Hom^i_{\catD(Y)}(\fmf{\cplx{K}}{X}{Y}(\cO_{Z_{x_1}}),\fmf{\cplx{K}}{X}{Y}(\cO_{x_2}))=0
$$
unless $x_1= x_2$ and $0\leq i\leq \dim X$.
\item There exists a closed point $x$ such that at least one of the following conditions is fulfilled:
\begin{enumerate} \item $\Hom^0_{\catD(Y)}(\fmf{\cplx{K}}{X}{Y}(\cO_X),
\fmf{\cplx{K}}{X}{Y}(\sk{x}))\simeq k$.
\item $\Hom^0_{\catD(Y)}(\fmf{\cplx{K}}{X}{Y}(\cO_{Z_x}),
\fmf{\cplx{K}}{X}{Y}(\sk{x}))\simeq k$ for any l.c.i.~zero cycle
$Z_x$ supported on $x$.
\item $1\leq \dim_k \Hom^0_{\catD(Y)}(\fmf{\cplx{K}}{X}{Y}(\cO_{Z_x}),
\fmf{\cplx{K}}{X}{Y}(\cO_{Z_x}))\leq l(\cO_{Z_x})$ for any
l.c.i.~zero cycle $Z_x$ supported on $x$, where $l(\cO_{Z_x})$ is
the length of $\cO_{Z_x}$ .
\end{enumerate}
\end{enumerate}
\end{defn}

Notice that now the objects $\bL j_{x}^\ast\cplx{K}$ are not required to be simple; hence we avoid the word simple in the denomination of the objects satisfying the above conditions and turn back to something closer to the original Bondal and Orlov way to describe them.

We prove now a variant of Theorem \ref{1:ffcritCM} which is valid in arbitrary characteristic.
Here, the requirement that $X$ is integral can be relaxed.

\begin{thm}\label{1:ffcritCMp}  Let $X$ and $Y$ be proper schemes over
 an algebraically closed field of arbitrary characteristic, and
 let $\cplx{K}$ be an object
in $\cdbc{X\times Y}$ of finite homological dimension over both
$X$ and $Y$. Assume also that $X$ is projective, Cohen-Macaulay,
equidimensional and connected. Then the functor
$\fmf{\cplx{K}}{X}{Y}\colon \cdbc{X}\to \cdbc{Y}$ is fully
faithful if and only if the kernel $\cplx{K}$ satisfy the Cohen-Macaulay orthonormality conditions over $X$ (Definition \ref{d:ortho}).
\end{thm}

\begin{proof}
The direct is immediate. For the converse we proceed as in the
proof of \ref{1:ffcritCM}. As there, using condition (1) of Definition \ref{d:ortho}, one sees that $\Phi$ has a
right adjoint  $H$ and that $H\circ \Phi\simeq\fmf{\cM}{X}{X}$,
where $\cM$ is a sheaf whose support is contained in  the diagonal
and $\pi_{1\ast}\cM$ is locally free. Since $X$ is connected, we can consider the rank
$r$ of
$\pi_{1\ast}\cM$, which is nonzero by
condition (2) of Definition \ref{d:ortho}; thus the support of $\cM$ is the diagonal. To
conclude, we have only to prove that $r=1$.

If $\cplx{K}$ satisfies (2.1) of Definition \ref{d:ortho}, then
$$\Hom^0_{\catD(X)}(\cO_X,\fmf{\cM}{X}{X}(\cO_x)) \simeq
\Hom^0_{\catD(Y)} (\fmf{\cplx{K}}{X}{Y}(\cO_X),
\fmf{\cplx{K}}{X}{Y}(\sk{x})) \simeq k.$$ Hence
$\fmf{\cM}{X}{X}(\cO_x)\simeq \cO_x$ and $r=1$.

If $\cplx{K}$ satisfies (2.2) of Definition \ref{d:ortho}, then
$$\Hom^0_{\catD(X)}(\cO_{Z_x},\fmf{\cM}{X}{X}(\cO_x)) \simeq
\Hom^0_{\catD(Y)} (\fmf{\cplx{K}}{X}{Y}(\cO_{Z_x}),
\fmf{\cplx{K}}{X}{Y}(\sk{x})) \simeq k$$ for any l.c.i.~zero cycle
$Z_x$. Hence $\fmf{\cM}{X}{X}(\cO_x)\simeq \cO_x$ and $r=1$.

Finally, assume that $\cplx{K}$ satisfies (2.3) of Definition \ref{d:ortho}, and let us prove that then condition (2.2) of Definition \ref{d:ortho} holds as well.

If $\widetilde{\cplx\cM} =\bR\dSHom{\cO_{X\times X}}(\cM,\pi_1^!\cO_X)$, proceeding as in the proof of
Theorem \ref{1:ffcritCM}, one has that
$\widetilde{\cplx\cM}$ is a sheaf $\widetilde\cM$ supported
topologically on the diagonal and that $\pi_{1\ast}\widetilde\cM$ is
locally free.  It follows that the functor $\fmf{\cM}{X}{X}$ has a
left adjoint $G$, defined as $G(\cplx
F)=\bR\pi_{1\ast}(\pi_2^\ast \cplx F\lotimes \widetilde \cM)$. This can be seen as follows: one has that
$\cM\simeq \bR\dSHom{\cO_{X\times X}}(\widetilde\cM, \pi_1^!\cO_X)$ by Proposition \ref{dsharp}, and
then $\pi_1^\ast \cplx G\lotimes \cM \simeq \bR\dSHom{\cO_{X\times X}}(\widetilde\cM, \pi_1^!\cplx G)$ by
Proposition \ref{dualGore};
thus
$$ \Hom (\cplx F, \fmf{\cM}{X}{X}(\cplx G)) \simeq \Hom (\pi_2^\ast\cplx F,
\bR\dSHom{\cO_{X\times X}}(\widetilde\cM, \pi_1^!\cplx G) ) \simeq
\Hom (G(\cplx F),\cplx G)\,,
$$
which proves that $G$ is a left adjoint to $\fmf{\cM}{X}{X}$. Thus, condition (2.2) is equivalent to
$\Hom^0_{\catD(X)}(G(\cO_{Z_x}),\cO_x)\simeq k$.

We know that $\cM$ is a sheaf topologically supported on the diagonal and
$\pi_{1\ast}\cM$ is locally free. Then,  if $\cF$ is a
sheaf, $\fmf{\cM}{X}{X}(\cF)$ is also a sheaf and the functor
$\cF\mapsto \phi(\cF)= \fmf{\cM}{X}{X}(\cF)$ is exact. One has that ${\calH}^i(G(\cF))=0$ for $i>0$ because $\widetilde \cM$ is topologically
supported on the diagonal. Hence,
$\Hom_{\catD(X)}(G(\cF_1),\cF_2) \simeq
\Hom_{\catD(X)}(G^0(\cF_1),\cF_2)$ for whatever sheaves $\cF_1$, $\cF_2$, where $G^0=\calH^0 \circ G$. This has two consequences: first, $G^0$ is a left adjoint to $\phi$; second, there are isomorphisms
$$
\Hom^0_{\catD(X)}(G(\cO_{Z_x}),\cO_x)\simeq
\Hom^0_{\catD(X)}(G^0(\cO_{Z_x}),\cO_x)\simeq
\Hom_{\cO_{Z_x}}(j_{Z_x}^\ast G^0(\cO_{Z_x}),\cO_x)\,,
$$
so that we are reduced to prove that $\Hom_{\cO_{Z_x}}(j_{Z_x}^\ast G^0(\cO_{Z_x}),\cO_x)\simeq  k$.
Then, it is enough to see that $j_{Z_x}^\ast G^0(\cO_{Z_x})\simeq \cO_{Z_x}$. 

Using the exactness of $\phi$, one proves by induction on  the length
$\ell(\cF)$ that the unit map $\cF\to \phi(\cF)$ is injective for any sheaf $\cF$ supported on $x$.  It follows
that the morphism $G^0(\cF)\to \cF$ is an epimorphism. To see that this is indeed the case, we have only to prove that for every sheaf $\cF'$ supported on the point $x$, the morphism $\Hom(\cF,\cF')\to \Hom(G^0(\cF),\cF')$ is injective. By the adjuntion formula, this is identified with the morphism $\Hom(\cF,\cF')\to  \Hom( \cF
,\phi(\cF'))$, which is injective because $\cF'\to \phi(\cF')$ is so.

In particular, 
 $\eta\colon G^0(\cO_{Z_x})\to \cO_{Z_x}$ is surjective, and
$\dim \Hom^0_{\catD(X)}(G^0(\cO_{Z_x}),\cO_{Z_x})\geq \ell (\cO_{Z_x})$. Thus, by condition (2.3) of Definition \ref{d:ortho},  $\dim \Hom^0_{\catD(X)}(G^0(\cO_{Z_x}),\cO_{Z_x}) = \ell(\cO_{Z_x})$. 
Since $\cO_{Z_x}$ is free, the exact sequence of
$\cO_{Z_x}$-modules
$$
0\to\cN\to j_{Z_x}^\ast G^0(\cO_{Z_x})\xrightarrow{j_{Z_x}^\ast(\eta)} \cO_{Z_x}\to 0
$$
splits, so that
$$
 0\to  \Hom_{\cO_{Z_x}}( \cO_{Z_x}  ,\cO_{Z_x})\to \Hom_{\cO_{Z_x}}( j_{Z_x}^\ast G^0(\cO_{Z_x}) ,\cO_{Z_x}) \to
\Hom_{\cO_{Z_x}}( \cN  ,\cO_{Z_x})\to 0 
$$
is an exact sequence. Moreover, 
$\Hom_{\cO_{Z_x}}( \cN ,\cO_{Z_x})=0$ because the two first terms have the same dimension. Let us see that this implies $\cN=0$. If $\cO_x\to\cO_{Z_x}$ is a nonzero, and then injective, morphism, we have $\Hom_{\cO_{Z_x}}(\cN,\cO_x)=0$ so that $\cN=0$ by Nakayama's lemma.
\end{proof}

\begin{rem} Even if the base field has characteristic zero, there is no obvious direct relationship between strongly simplicity and the Cohen-Macaulay orthonormality conditions for a kernel $\cplx K$ in $\cdbc{X\times Y}$, despite the fact that they are equivalent due to Theorems \ref{1:ffcritCM} and \ref{1:ffcritCMp}. Then, Theorem \ref{1:ffcritCMp} is actually a new characterisation of which kernels induce a fully faithful integral functor. 

If we wish to compare directly strongly simplicity with the Cohen-Macaulay orthonormality conditions, we see that the first condition in both definitions is the same; however, the second condition of strongly simplicity is a property that has to be satisfied at \emph{every} closed point for \emph{one} cycle supported on the point, whereas two of the forms of the second condition of the Cohen-Macaulay orthonormality conditions refer to a property which has to be satisfied for \emph{all} zero cycles supported on only \emph{one}  closed point. 
\end{rem}
\subsection{A criterion in the relative setting\label{ellfibss}}

In the relative situation the notions of strongly simple object  and of an object satisfying the  Cohen-Maculay orthonormality conditions, are
the following.

\begin{defn} Assume that $X\to S$ is
Cohen-Macaulay. An object $\cplx{K}\in \cdbc{X\times_SY}$ is
\emph{relatively strongly simple (resp.~satisfies the relative Cohen-Maculay orthonormality conditions)} over $X$ if $\cplx{K}_s$ is
strongly simple (resp.~satisfies the Cohen-Maculay orthonormality conditions) over $X_s$  for every closed point $s\in S$.
\end{defn}

As a corollary of  Proposition \ref{p:relative} and Theorems
\ref{1:ffcritCM} and \ref{1:ffcritCMp}, we obtain the following
result.

\begin{thm} \label{t:relative}
Let $X\to S$ and $Y\to S$ be proper and flat
 morphisms. Assume also that $ X\to S$ is locally projective with
 Cohen-Macaulay, equidimensional, and connected fibers.  Let $\cplx{K}$ be an object in
$\cdbc{X\times_SY}$ of finite homological dimension over both $X$
and $Y$. \begin{enumerate} \item Assume that we are in
characteristic 0 and $X \to S$ has integral fibers. The relative
integral functor $\fmf{\cplx{K}}{X}{Y}\colon \cdbc{X}\to \cdbc{Y}$
is fully faithful  if and only if $\cplx{K}$ is relatively
strongly simple over $X$. \item In arbitrary characteristic,
$\fmf{\cplx{K}}{X}{Y}\colon \cdbc{X}\to \cdbc{Y}$ is fully
faithful if and only if $\cplx{K}$ satisfies the relative Cohen-Macaulay orthonormality conditions over $X$.
\end{enumerate}
\end{thm}

\subsubsection{Application to genus one fibrations in arbitrary characteristic}

We now apply Theorem \ref{t:relative} to give an alternative proof for
Proposition \ref{p:poincequiv} without using spherical objects.

With the same notation than in the Subsection \ref{ss:autoequiv}
and taking into account  the symmetry of  $\cI_\Delta$,  to
conclude that the functor $\Phi=\fmf{\cI_\Delta}{X}{X}$ is an
auto-equivalence of $\cdbc{X}$, it is enough to prove that $\cI_\Delta$ satisfies the relative
Cohen-Macaulay orthonormality conditions over  the first factor.
We fix a closed point $s\in S$ and consider two points $x$ and
$\bar x$  in the fiber $X_s$.

Let $Z_x\hookrightarrow X_s$ be a l.c.i.~zero cycle supported on $x$ of length
$\ell$ defined by an ideal  $\cI_{Z_x}$; then  $\cI_{Z_x}$ is an invertible
sheaf of $\cO_{X_s}$-modules. We denote by $\pi\colon Z_x\times X_s\to X_s$ the second projection and by $\cJ_{Z_x}$ the push-forward by $\pi$ of the ideal of the graph of
$Z_x\hookrightarrow X_s$. We have a commutative diagram of exact rows
\begin{equation} \label{e:diagram}
\xymatrix{0 \ar[r] & \cI_{Z_x} \ar[r]\ar@{^(->}[d] & \cO_{X_s}  \ar[r]\ar@{^(->}[d] & \cO_{Z_x}  \ar[r]\ar@{=}[d] & 0 \\
0\ar[r] & \cJ_{Z_x} \ar[r] &\pi_\ast\cO_{Z_x\times X_s} \ar[r] &  \cO_{Z_x} \ar[r] & 0
}
\end{equation}
Since $\pi_\ast\cO_{Z_x\times X_s}$ is a free $\cO_{X_s}$-module of rank $\ell$, and $\pi$ has a section, one easily sees that the quotient sheaf $\pi_\ast\cO_{Z_x\times X_s} / \cO_{X_s}$ is free of rank $\ell-1$ and we have an exact sequence
$$
0\to \cI_{Z_x}\to \cJ_{Z_x} \to \cO_{X_s}^{\oplus(\ell-1)}\to 0\,,
$$
which proves that $\cJ_{Z_x}$ is locally free of rank $\ell$.

One has that $\Phi_s(\cO_{Z_x})\simeq \cJ_{Z_x}$ and $\Phi_s(\cO_x)\simeq \cI_x$, where $\cI_x$ is the ideal of the point $x$;  we have then to compute the groups  $\Hom^i_{\catD(X_s)}(\cJ_{Z_x}, \cI_{\bar x})$. Since $\cJ_{Z_x}$ is locally free, we know that
$\Hom^i_{\catD(X_s)}(\cJ_{Z_x}, \cI_{\bar x})=0$ unless  $i=0,1$, so we have to worry only about the cases $i=0$ and $i=1$.

We have that
$$
\Hom^0(\pi_\ast\cO_{Z_x\times X_s}, \cI_{\bar x})\simeq \Hom^0(\cO_{X_s}^{\oplus\ell}, \cI_{\bar x})\simeq H^0(X_s,\cI_{\bar x})^{\oplus\ell}=0\,.
$$
Then, from the bottom row of diagram \eqref{e:diagram} we obtain the following diagram
{\footnotesize
\begin{equation} \label{e:diagram2}
\xymatrix@C=8pt{ 0\ar[r] &  \Hom^0(\cJ_{Z_x},\cI_{\bar x}) \ar[r] \ar@{^(->}[d] &  \Hom^1(\cO_{Z_x}, \cI_{\bar x}) \ar[r] \ar[d]& \Hom^1(\pi_\ast\cO_{Z_x\times X_s},
\cI_{\bar x})\ar[r] \ar[d]^h & \Hom^1(\cJ_{Z_x},\cI_{\bar x})\ar[r] \ar[d] & 0\\
\dots \ar[r]&  \Hom^0( \cJ_{Z_x},\cO_{X_s}) \ar[r]  &  \Hom^1(\cO_{Z_x}, \cO_{X_s}) \ar[r]^(.43)g& \Hom^1(\pi_\ast\cO_{Z_x\times X_s}, \cO_{X_s})\ar[r] & \Hom^1( \cJ_{Z_x} ,\cO_{X_s})\ar[r] & \dots
}
\end{equation}
}

\medskip\noindent
a) \emph{The morphisms $h$ and $g$ are isomorphisms.}
\medskip

The fact that $h$  is an isomorphism follows from the formula $\pi_\ast\cO_{Z_x\times X_s}\simeq
\cO_{X_s}^{\oplus\ell}$.
To prove that $g$ is also an isomorphism, we first notice that $g$ is the morphism obtained by applying
the functor $\Hom^1(-,\cO_{X_s})$ to the projection $\varpi\colon \pi_\ast\cO_{Z_x\times X_s}\to \cO_{Z_x}$. By duality, $\Hom^1(-,\cO_{X_s})\simeq H^0(X_s,-)^\ast$ and one finishes because $\varpi\colon \pi_\ast\cO_{Z_x\times X_s}\to
\cO_{Z_x}$  induces a isomorphism between the corresponding spaces of global sections.

\medskip\noindent
b) \emph{Condition (1) of Definition \ref{d:ortho}}
\medskip

It suffices to see that
$\Hom^1(\cO_{Z_x}, \cI_{\bar x})\to \Hom^1(\pi_\ast\cO_{Z_x\times X_s}, \cI_{\bar x})$ is an
isomorphism when $ x\neq \bar x$.  In this case,  the second vertical arrow of diagram \eqref{e:diagram2} is an isomorphism, so that the above morphism is identified with $g$, which is an isomorphism.

\medskip\noindent
c) \emph{Condition (2.2) of Definition \ref{d:ortho}}
\medskip

We have to prove that
$\Hom^0(\cJ_{Z_x},\cI_x)\simeq k$.  From diagram \eqref{e:diagram2} $\Hom^0(\cJ_{Z_x},\cI_x)$ is isomorphic to the kernel of the second vertical arrow $\Hom^1(\cO_{Z_x},\cI_x) \to \Hom^1(\cO_{Z_x},\cO_{X_s})$, which is isomorphic to $\Hom^0(\cO_{Z_x},\cO_x)\simeq k$.

Notice that the proof of this result does not use spanning classes
and is valid in any characteristic.


\section{Fourier-Mukai partners}\label{intpartnersSec}

When two proper schemes $X$ and $Y$ have equivalent derived
categories $\cdbc{X}\simeq \cdbc{Y}$, they are called
\emph{$D$-equivalent}. When the equivalence is given by an
integral functor  (such functors are called  Fourier-Mukai
functors), we have the following more restrictive notion.

\begin{defn} \label{d:intpartners}
Two proper schemes $X$ and $Y$ are \emph{Fourier-Mukai partners} if there is a Fourier-Mukai functor
$$
\fmf{\cplx K}{X}{Y}\colon \cdbc{X}\simeq \cdbc{Y}\,.
$$
\end{defn}

Due to Orlov's representation theorem  \cite{Or97}, if $X$ and $Y$
are smooth and projective then they are $D$-equivalent if and only
if they are Fourier-Mukai partners.  Since the validity of Orlov's
theorem for singular varieties is still unknown, we shall adopt
the (hopefully provisional) notion of Fourier-Mukai partners.

In the rest of the section $X$ and $Y$ are projective and $\cplx K$ is a kernel
in $\cdbc{X\times Y}$  such that $\fmf{\cplx K}{X}{Y}$ is an equivalence $\cdbc{X}\simeq \cdbc{Y}$.

\begin{prop} \label{dualsiso}
There is a natural isomorphism
$$
\bR \dSHom{\cO_{X\times Y}}(\cplx K, \pi_X^!\cO_X)\simeq \bR \dSHom{\cO_{X\times Y}}(\cplx K, \pi_Y^!\cO_Y)\,.
$$
\end{prop}
\begin{proof} By Proposition \ref{p:equivfhd} the right adjoint to
$\fmf{\cplx{K}}{X}{Y}$ is $\fmf{\bR\dSHom{\cO_{X\times
Y}}(\cplx{K}, \pi_Y^!\cO_Y) }{Y}{X}$ and the right adjoint to the latter functor
 is  $\fmf{\cplx L}{X}{Y}$ where
 $$
 \cplx L=
\bR\dSHom{\cO_{X\times Y}}(\bR\dSHom{\cO_{X\times Y}}(\cplx{K},
\pi_Y^!\cO_Y), \pi_X^!\cO_X)\,.
$$
Since the left and right adjoint of
an equivalence are naturally isomorphic, one has  that $\fmf{\cplx K}{X}{Y}\simeq
\fmf{\cplx L}{X}{Y}$.

Moreover, for any scheme $T$ we can consider $X_T=X\times T$ and
$Y_T=Y\times T$. The kernel $\cplx K$ gives rise to a relative
kernel $\cplx{K}_T =\pi_{X\times Y}^\ast \cplx K \in
\cdbc{X_T\times_T Y_T}$, where $\pi_{X\times Y}$ is the projection
$X_T\times_T Y_T\simeq X\times Y\times T \to X\times Y$. This
relative kernel is of finite homological dimension over both $X_T$
and $Y_T$, and the relative integral functor
$\fmf{\cplx{K}_T}{X_T}{Y_T} \colon \cdbc{X_T}\to\cdbc{Y_T}$ is an
equivalence by Proposition \ref{p:relative}. Arguing as above,  we
get $\fmf{\cplx{K}_T}{X_T}{Y_T}\simeq \fmf{\cplx L_T}{X_T}{Y_T}$.
If we take $T=X$ and apply the above isomorphism to the sheaf
$\cO_\Delta$ of the diagonal, we obtain an isomorphism
$$
\cplx K\simeq \bR\dSHom{\cO_{X\times Y}}(\bR\dSHom{\cO_{X\times Y}}
(\cplx{K}, \pi_Y^!\cO_Y), \pi_X^!\cO_X)\,.
$$
Since $\bR\dSHom{X\times Y}(\cplx{K}, \pi_Y^!\cO_Y)$ is of finite
homological dimension over $X$, we conclude by Proposition \ref{dsharp}.
\end{proof}

\begin{rem} When $X$ and $Y$ are smooth, Proposition \ref{dualsiso} is equivalent to the fact
that any equivalence commutes with the Serre functors. Then,
Proposition \ref{dualsiso} can be  considered as a generalisation
to that property for arbitrary singular schemes.
\end{rem}

Our next aim is to prove that any integral Fourier-Mukai partner
of a projective Cohen-Macaulay  (resp.~Gorenstein) scheme $X$ is
also Cohen-Macaulay  (resp.~Gorenstein).



\begin{thm} \label{t:intpartnersCM}
Let $X$ be a projective equidimensional Cohen-Macaulay scheme and
$Y$ a projective Fourier-Mukai partner  of $X$. Then one has
\begin{enumerate}
\item If $Y$ is reduced, then $Y$ is equidimensional of dimension $m=\dim X$.
\item If $Y$ is equidimensional and $\dim Y=\dim X$, then $Y$ is Cohen-Macaulay. Moreover, if $X$ is
Gorenstein, then $Y$ is Gorenstein as well.
\end{enumerate}
\end{thm}
\begin{proof} By Proposition \ref{dualsiso}
the integral functors
\begin{align*}
\Psi_1=\fmf{\bR \dSHom{\cO_{X\times Y}}(\cplx K, \pi_X^!\cO_X)}{X}{Y}\colon \cdbc{X} & \to \cdbc{Y} \\
\Psi_2=\fmf{\bR \dSHom{\cO_{X\times Y}}(\cplx K,
\pi_Y^!\cO_Y)}{X}{Y}\colon \cdbc{X} & \to \cdbc{Y}\,,
\end{align*}
are naturally isomorphic. Hence, for any l.c.i zero cycle $Z_x$ supported on a
closed point $x\in X$ and any closed point $y\in Y$ we have
$$
\bR\dSHom{\cO_Y}(\cO_y, \Psi_1(\cO_{Z_x}) ) \simeq
\bR\dSHom{\cO_Y}(\cO_y, \Psi_2(\cO_{Z_x}) )\,.
$$
This gives rise to the formula:
\begin{equation} \label{keyiso}
\bR \dSHom{\cO_Y}(\Phi_{Z_x}(\Dcplx{Z_x}),\cO_y)\simeq \bR
\dSHom{\cO_Y}(\Phi_{Z_x} (\cO_{Z_x}),\bR
\dSHom{\cO_Y}(\cO_y,\cO_Y))[m]\,,
\end{equation}
 where $\Phi_{Z_x}\colon \cdbc{Z_x} \to
\cdbc{Y}$ is the integral functor of kernel $\bcplx{K}= \bL
j_{Z_x}^\ast\cplx K\in \cdbc{Z_x\times Y}$  (see Lemmas \ref{Psi1} and \ref{Psi2} for details).

Since $X$ is Cohen-Macaulay, every l.c.i.~cycle $Z_x$ is
Cohen-Macaulay as well, and its dualizing complex $\Dcplx{Z_x}$ is
a single sheaf $\omega_{Z_x}$.

As $\Phi$ is an equivalence of
categories, it follows from Lemma \ref{l:spanning} that the
objects $\Phi_{Z_x}(\cO_{Z_x})\simeq\Phi(\cO_{Z_x})$ for all
l.c.i.~cycles $Z_x$ form a spanning class for $\cdbc{Y}$. Then,
if we fix a closed point  $y\in Y$,
there is a l.c.i.~cycle $Z_x$ such that $y\in
\supp(\Phi_{Z_x}(\cO_{Z_x}))$.

Let us denote by  $\phi\colon Y'=\Spec \cO_{Y,y}\to Y$ the
natural flat morphism.
The composition  $\Phi'_{Z_x}=\phi^\ast\circ \Phi_{Z_x}$ is
an integral functor with kernel
$\bcplx{K'}=(1\times\phi)^\ast \bcplx{K}\in \cdbc{Z_x\times Y'}$.
If we apply
$\phi^\ast $ to \eqref{keyiso} we get the analogous formula
\begin{equation} \label{keyiso2}
\bR \dSHom{\cO_{Y'}}(\Phi'_{Z_x}(\omega_{Z_x}),\cO_y)\simeq \bR
\dSHom{\cO_{Y'}}(\Phi'_{Z_x} (\cO_{Z_x}),\bR
\dSHom{\cO_{Y'}}(\cO_y,\cO_{Y'}))[m]\,,
\end{equation}

After shifting degrees if necessary, we can  assume that
$\calH^0(\bcplx {K'})\neq 0$ and $\calH^i(\bcplx {K'})= 0$ for
$i>0$. Then $\calH^0(\Phi'_{Z_x}(\cO_{Z_x}))\neq 0$ and
$\calH^i(\Phi'_{Z_x}(\cO_{Z_x}))= 0$ for $i>0$. For every nonzero finite
$\cO_{Z_x}$-module $\cF$, one has  that $\calH^0
(\pi_{Z_x}^\ast \cF\lotimes \bcplx {K'})\simeq \pi_{Z_x}^\ast \cF\otimes
\calH^0(\bcplx {K'}) \neq 0$ and that $\calH^i(\pi_{Z_x}^\ast \cF\lotimes
\bcplx {K'})=0$ for $i>0$, where $\pi_{Z_x}\colon Z_x\times
Y'\to Z_x$ is the projection. Thus, $\calH^0(\Phi'_{Z_x}(\cF))\neq
0$ and $\calH^i(\Phi'_{Z_x}(\cF))= 0$ for $i>0$.

It follows that $\bR\dSHom{\cO_{Y'}}(\Phi'_{Z_x} (\omega_{Z_x}),
\cO_y)$ has no negative cohomology  sheaves and that its $0$-th
cohomology sheaf is nonzero.
By Equation \eqref{keyiso2}, the same happens for the object
$\cplx Q=\bR\dSHom{\cO_{Y'}}(\Phi'_{Z_x} (\cO_{Z_x}),
\bR\dSHom{\cO_{Y'}}(\cO_y,\cO_{Y'}) )[m]$, that is,
\begin{equation}
\label{e:keyiso2} \calH^i(\cplx Q)=0 \text{\ for $i<0$, and\ }
\calH^0(\cplx Q)\neq 0\,.
\end{equation}

\medskip

(1) Assume that $Y$ is reduced and we choose $y$ to be a smooth
point of an irreducible  component $Y_0$ of $Y$. Then
$\bR\dSHom{\cO_{Y'}}(\cO_y,\cO_{Y'})\simeq \cO_y[-n]$, where
$n=\dim Y_0$, and we have  $\cplx Q\simeq
\bR\dSHom{\cO_{Y'}}(\Phi'_{Z_x}(\cO_{Z_x}), \cO_y )[m-n]$. Since
both $\cplx Q$ and
$\bR\dSHom{\cO_{Y'}}(\Phi'_{Z_x}(\cO_{Z_x}),\cO_y)$ have no
negative cohomology sheaves and nonzero 0-th cohomology sheaf, one
must have $m-n=0$. Thus all the irreducible components of $Y$ have
dimension $m$.

\medskip

(2)  Let $y\in Y$ be a closed point  and let  $j_0$ be the first
index $j$ with $ \SExt^j_{\cO_{Y'}}(\cO_y,\cO_{Y'})\neq 0$, that
is, $j_0$ is the depth of the local ring $\cO_{Y,y}$. Then
$$
\calH^{j_0} [\cplx Q[-m]) \simeq \SHom_{\cO_{Y'}} (\calH^0
(\Phi'_{Z_x} (\cO_{Z_x})), \SExt^{j_0}_{\cO_{Y'}}(\cO_y,\cO_{Y'}))
\,,
$$
  which
is not zero because $\SExt^{j_0}_{\cO_{Y'}}(\cO_y,\cO_{Y'})$ is a
nonzero $\cO_y$-vector space. Hence $j_0\geq m=\dim Y$
(cf.~Equation \eqref{e:keyiso2}),  and $Y$ is Cohen-Macaulay.

Assume now that $X$ is Gorenstein. Then $\omega_{Z_x}\simeq
\cO_{Z_x}$ and we have
$$
\bR \dSHom{\cO_{Y'}}(\Phi'_{Z_x}(\cO_{Z_x}),\cO_y)\simeq\bR
\dSHom{\cO_{Y'}}(\Phi'_{Z_x}(\cO_{Z_x}),\bR
\dSHom{\cO_{Y'}}(\cO_y,\cO_{Y'}))[m]\,,
$$
by Equation \eqref{keyiso2}. Since
$\SExt^j_{\cO_{Y'}}(\cO_y,\cO_{Y'})=0$ for $j<m$, we deduce that
$$
\Hom_{Y'}(\calH^0(\Phi'_{Z_x}(\cO_{Z_x})),\cO_y)\simeq
\Hom_{Y'}(\calH^0(\Phi'_{Z_x}(\cO_{Z_x})),
\SExt^m_{\cO_{Y'}}(\cO_y,\cO_{Y'}))\,,
$$
and then,
$$
\Hom_{\cO_y}(\calH^0(\Phi'_{Z_x}(\cO_{Z_x}))\otimes
\cO_y,\cO_y)\simeq
\Hom_{\cO_y}(\calH^0(\Phi'_{Z_x}(\cO_{Z_x}))\otimes \cO_y,
\SExt^m_{\cO_{Y'}}(\cO_y,\cO_{Y'}))\,.
$$
Thus, $\dim \SExt^m_{\cO_{Y'}}(\cO_y,\cO_{Y'})=1$ and $Y$ is
Gorenstein.
\end{proof}

We now prove the auxiliary Lemmas used in the proof of Theorem \ref{t:intpartnersCM}.

\begin{lem}\label{Psi1}  $\bR \dSHom{\cO_Y}(\cO_y, \Psi_1(\cO_{Z_x}) ) \simeq \bR
\dSHom{\cO_Y}(\Phi_{Z_x}(\Dcplx{Z_x}),\cO_y)$.
\end{lem}

\begin{proof} We first compute $\Psi_1(\cO_{Z_x})$. By Corollary \ref{dsharpequiv}
and Lemma \ref{sharpfmf}
$$
\Psi_1(\cO_{Z_x}) \simeq \Psi_1(\cO_{Z_x}^{\sharp\sharp})\simeq
(\Phi(\cO_{Z_x}^\sharp))^\sharp\simeq
(\Phi_{Z_x}(\Dcplx{Z_x}))^\sharp\,. $$ where the last isomorphism
is due to the isomorphism $\cO_{Z_x}^\sharp \simeq j_{Z_x\ast}
\Dcplx{Z_x}$, which follows by relative duality for the closed
immersion $j_{Z_x}\colon Z_x\hookrightarrow X$. Considering
Equation \eqref{tens1}, we obtain:
$$\aligned \bR\dSHom{\cO_Y}(\cO_y, \Psi_1(\cO_{Z_x}) ) &\simeq
\bR\dSHom{\cO_Y}(\cO_y, (\Phi_{Z_x}(\Dcplx{Z_x}))^\sharp)\\
&\simeq \bR \dSHom{\cO_Y}(\cO_y\lotimes
\Phi_{Z_x}(\Dcplx{Z_x}),\Dcplx{Y}) \simeq \bR
\dSHom{\cO_Y}(\Phi_{Z_x}(\Dcplx{Z_x}),\cO_y^\sharp)\endaligned $$
Relative duality for the closed immersion $\{y\}\hookrightarrow Y$
gives $\cO_y^\sharp\simeq \cO_y$, which finishes the proof.
\end{proof}

\begin{lem}\label{Psi2} One has $$\bR \dSHom{\cO_Y}(\cO_y, \Psi_2(\cO_{Z_x}) ) \simeq \bR
\dSHom{\cO_Y}(\Phi_{Z_x} (\cO_{Z_x}),\bR
\dSHom{\cO_Y}(\cO_y,\cO_Y))[m]\,.
$$
\end{lem}

\begin{proof} We first compute $\Psi_2(\cO_{Z_x})$.  Let us denote $\pi_{Z_x}$ and $\bar\pi_Y$
 the projections of $Z_x\times Y$ onto $Z_x$ and $Y$ respectively. Using that $\cO_{Z_x\times Y}$ is of
finite homological dimension as a module over $\cO_{X\times Y}$
and taking into account Equation \eqref{tens2}, we have that
$$
\Psi_2(\cO_{Z_x}) \simeq \bar\pi_{Y\ast} \bR \dSHom{\cO_{Z_x\times
Y}}(\bcplx K, \bL j_{Z_x}^\ast\pi_Y^!\cO_Y)\,.
$$
The term $\bL j_{Z_x}^\ast\pi_Y^!\cO_Y$ is computed as $\bL
j_{Z_x}^\ast\pi_{X}^\ast \Dcplx{X}\simeq  \pi_{Z_x}^\ast \bL
j_{Z_x}^\ast \Dcplx{X}$. Furthermore, relative duality for the
regular immersion $j_{Z_x}\colon Z_x \hookrightarrow X$ gives $
j_{Z_x}^! \cO_X\simeq \cO_{Z_x}[-m]$, where $m=\dim X$ and
$\Dcplx{Z_x}\simeq j_{Z_x}^! \cO_X \lotimes \bL j_{Z_x}^\ast
\Dcplx{X}$. Thus, $\pi_{Z_x}^\ast \bL j_{Z_x}^\ast \Dcplx{X}\simeq
\pi_{Z_x}^\ast \Dcplx{Z_x}[m]\simeq \bar\pi_Y^! \cO_Y[m]$, so that
$$
\begin{aligned}
\Psi_2(\cO_{Z_x})&\simeq \bar \pi_{Y\ast} \bR \dSHom{\cO_{Z_x\times Y}}(\bcplx K,  \bar\pi_Y^!\cO_Y[m])\\
& \simeq
 \bR \dSHom{\cO_Y}(\Phi_{Z_x}(\cO_{Z_x}),\cO_Y[m])
 \,.
\end{aligned}$$ Then $$\bR \dSHom{\cO_Y}(\cO_y, \Psi_2(\cO_{Z_x})
)\simeq \bR \dSHom{\cO_Y}(\cO_y, \bR
\dSHom{\cO_Y}(\Phi_{Z_x}(\cO_{Z_x}),\cO_Y[m])$$ and one concludes
by Equation \eqref{tens1}.

\end{proof}



\end{document}